\documentclass{amsart}
\usepackage{latexsym,amsmath,amssymb,amsfonts,amscd,graphics}
\usepackage[mathscr]{eucal}
\usepackage[all]{xypic}

\numberwithin{equation}{section}
\theoremstyle{plain}
\newtheorem{theorem}[equation]{Theorem}

\newtheorem{lemma}[equation]{Lemma}
\newtheorem{proposition}[equation]{Proposition}
\newtheorem{corollary}[equation]{Corollary}
\newtheorem{cordef}[equation]{Corollary/Definition}

\theoremstyle{remark}
\newtheorem{remark}[equation]{Remark}

\theoremstyle{definition}
\newtheorem{definition}[equation]{Definition}
\newtheorem{notation}[equation]{Notation}
\newtheorem{convention}[equation]{Convention}
\newtheorem{example}[equation]{Example}

\newcommand{\om}{\omega}
\newcommand{\omb}{\overline{\omega}}
\newcommand{\bP}{\mathbb{P}}

\newcommand{\bR}{\mathbb{R}}

\newcommand{\bQ}{\mathbb{Q}}
\newcommand{\bZ}{\mathbb{Z}}

\newcommand{\bC}{\mathbb{C}}

\newcommand{\bL}{\mathbb{L}}

\newcommand{\sH}{\mathrm{H}}

\newcommand{\calE}{\mathcal{E}}

\newcommand{\calO}{\mathcal{O}}

\newcommand{\calX}{\mathcal{X}}
\newcommand{\calB}{\mathcal{B}}

\newcommand{\calD}{\mathcal{D}}

\newcommand{\Span}{\mathrm{span}}

\newcommand{\Sat}{\mathrm{Sat}}

\newcommand{\s}{\mathrm{sp}}

\newcommand{\Gr}{\mathrm{Gr}}

\newcommand{\rank}{\mathrm{rank}}

\newcommand{\Pic}{\mathrm{Pic}}

\newcommand{\re}{\mathrm{Re}}
\newcommand{\im}{\mathrm{Im}}

\author{Radu Laza}
\title[]{Triangulations of the sphere and degenerations of  $K3$ surfaces} 
\address{University of Michigan \\
3863 East Hall \\
Ann Arbor, MI 48109}
\email{rlaza@umich.edu}
\begin{document}
\bibliographystyle{amsplain}

%%%%%%%%%%%%%%%%%%%%%%%%%%%%%%%%%%%%%%%%%%%%%%%%%%%%%%%%%%%%%%
\begin{abstract}
W. Thurston \cite{thurston} proved that to a triangulation of $S^2$ of non-negative combinatorial curvature, one can associate an element in a certain lattice over the Eisenstein integers such that its orbit is a complete invariant of the triangulation. In this paper, we show that this association can be obtained naturally by using Type III degenerations of $K3$ surfaces.  From this perspective, Thurston's result  can be interpreted as a hint towards the construction of a geometrically meaningful compactification for the moduli space of polarized $K3$ surfaces.
\end{abstract}
\maketitle

%%%%%%%%%%%%%%%%%%%%%%%%%%%%%%%%%%%%%%%%%%%%%%%%%%%%%%%%%%%%%%
%%% S0
\section*{Introduction}
It is well known that in order to find a good geometric compactification for the moduli space of $K3$ surfaces it is important  to understand the combinatorics of the triangulations of $S^2$ (see \cite[pg. 22--24]{friedmanmorrison}). In fact,  it is reasonable to expect that if a good compactification exists then the set of all triangulations of $S^2$ has an arithmetic structure related to the arithmetic of the $K3$ lattice. The latter is indeed true, as shown by W. Thurston \cite{thurston}. Specifically, Theorem 0.1 of \cite{thurston}  says that the triangulations of non-negative combinatorial curvature of $S^2$ (i.e. each vertex of the triangulation has degree at most $6$) are parameterized by (the orbits of) the points of positive norm in a lattice $M^\calE$ over the Eisenstein integers $\calE$. In this paper, we show that this result is natural in the context of degenerations of $K3$ surfaces. The purpose of doing this is twofold:  we clarify certain aspects of \cite{thurston} (e.g. we identify the lattice $M^\calE$), and, at the same time shed some light on the combinatorics of degenerations of  $K3$ surfaces.% (with an eye on the compactification problem).

\smallskip

The triangulations of $S^2$ occur in the study of the moduli space of $K3$ surfaces as the dual the graph of Type III degenerations of $K3$ surfaces. Specifically, as a consequence of the Kulikov--Persson--Pinkham theorem,  any $1$-parameter degeneration $\calX^*\to \Delta^*$ can be filled-in (possibly after a base change) to a semi-stable family $\calX\to \Delta$ with $K_\calX$  trivial. Depending on the index of nilpotence of the monodromy of the degeneration, the central fibers are classified  in three classes: Type I, II, and III. The Type I and II cases are relatively easily understood  (see \cite{friedmanmorrison}) and we will not discuss them here. In the Type III case, the central fiber $X_0$ is a normal crossing variety, whose dual graph is a triangulation of $S^2$. The normalizations of the components $V_i$ of $X_0$ are rational surfaces, meeting the other components along anticanonical cycles $D_i=\cup D_{ij}$ of rational curves (where $D_{ij}=V_i\cap V_j$ are the double curves). 
The isotopy class of the triangulation of $S^2$ given by the dual graph of $X_0$ and the self-intersection numbers $D^2_{ij}$ are important discrete invariants associated to a Type III $K3$ surface. By results of Friedman \cite{friedmansmooth}, any triangulation of $S^2$ and any combinatorially allowable assignment of the self-intersection numbers $D_{ij}^2$ can be realized for some Type III degeneration. 
  
 \smallskip

In this paper, we restrict our attention to Type III degenerations of $K3$ surfaces in minus-one-form (see Miranda--Morrison \cite{mm}), i.e. degenerations for which each double curve $D_{ij}$ has self-intersection $-1$. The point being that, due to the symmetry, the combinatorial data of a degeneration $X_0$ in  minus-one-form is completely determined by the isotopy class of the associated triangulation $T$ of $S^2$. Since for a degeneration in  minus-one-form, the associated triangulation is of non-negative combinatorial curvature, and each such triangulation can be realized, we conclude that the classification of non-negatively curved triangulations is equivalent to the classification of locally trivial deformation classes of  Type III degenerations  in minus-one-form. The advantage of considering degenerate $K3$ surfaces is that we can linearize the classification problem by taking the cohomology. 

\smallskip

The structure of the cohomology of a Type III $K3$ surface $X_0$ and the arithmetic involved in a Type III degeneration are quite well understood, esp. by work of Friedman--Scattone \cite{friedmanscattone}. Specifically, what is relevant for us is that to a type III $K3$ surface $X_0$ one can naturally associate a lattice $\overline{L}$ of signature $(1,18)$, representing  the cohomology classes that behave well with respect to a smoothing.  Furthermore, due to the existence of a smoothing to a $K3$ surface, $\overline{L}$ embeds into the lattice $M=E_8^{\oplus 2}\oplus U^{\oplus 2}$ with orthogonal complement spanned (up to some primitivity index $k$) by an element $\delta\in M$ of norm $t$, where $t$ is the number of triple points of $X_0$.  In the minus-one-case, the isometry class of the lattice $\overline{L}$ can be regarded as an intrinsic invariant associated to the triangulation given by the dual graph. Unfortunately, as shown by \cite[\S1]{friedmanscattone}, the only information contained in $\overline{L}$ is $t$ (and the index $k$), the number of triple points, or equivalently the number of triangles in the triangulation. Our main result in this paper (theorems \ref{mainthm1} and \ref{mainthm2}) is that we can naturally enrich the structure of $\overline{L}$ such that we obtain the  same arithmetic structure as Thurston \cite{thurston}, and thus recover the triangulation.

\smallskip

Our first step is to note that a Type III $K3$ surface $X_0$ in minus-one-form is canonically polarized by  the sum $h=\sum D_i$ of the anticanonical divisors of the components $V_i$. The polarization class $h$ can be regarded as an element of the lattice $\overline{L}$. By considering the primitive part $P=\langle h\rangle^\perp_{\overline L}$ of $\overline{L}$, we obtain a 
tower of lattice embeddings $P\hookrightarrow \overline{L}\hookrightarrow M$ that depends only on the  triangulation of $S^2$ given by the dual graph.  What is gained by considering the polarization is that $P$ is a negative definite lattice, containing  significantly more arithmetic information than $\overline{L}$. For example, in some cases, it is possible to distinguish between $2$ triangulations with the same number of triangles (see Ex. \ref{exdeg}).  

\smallskip

The second step of our construction is to show  that, due to the symmetry of the minus-one-form, the lattice $P$ comes equipped with a fixed-point-free isometry $\rho$ of order $3$. Even more, this isometry can be extended in a  compatible way to $M$.  We obtain that  $P$ and $M$ are lattices over the Eisenstein integers and the embedding $P\hookrightarrow M$ is an embedding of Eisenstein lattices.  By construction $h$ belongs to the orthogonal complement of $P$ in $M$. In order to make the construction compatible with the Eisenstein lattice structure, we replace $h$ by a naturally defined rank $2$ (or rank $1$ over $\calE$)  lattice $\Delta\subseteq P^\perp_M$ with $h\in \Delta$ and $\Delta\cong A_2(-\frac{t}{2})$. To emphasis the extra Eisenstein lattice structure we denote $P$, $M$, and $\Delta$ by $P^\calE$, $M^\calE$, and $\Delta^\calE$ respectively. It is then easily seen that entire arithmetic structure that we associated to a Type III $K3$ surface in minus-one-form can be conveniently encoded in the choice of a generator $\delta^\calE\in M^\calE$ for the rank $1$ sublattice $\Delta^\calE$. We note that the lattice $M^\calE$ is a standard Eisenstein lattice, with underlying $\bZ$-lattice $E_8^{\oplus 2}\oplus U^{\oplus 2}$, that occurs quite frequently in algebraic geometry, see esp. Allcock \cite{allcockaut}.

\smallskip

To conclude, to a triangulation of $S^2$ of non-negative combinatorial curvature, we associate (essentially canonically)  a polarized Type III $K3$ surface $(X_0,h)$ in minus-one-form.  Then, by considering the cohomology of $X_0$, we obtain an element $\delta^\calE$ in a lattice $M^\calE$ over the Eisenstein integers.  By construction, the orbit of $\delta^\calE$ in $M^\calE$ (w.r.t. the group of isometries $\Gamma^\calE$) is an arithmetic invariant of the triangulation. On the other hand, Thurston \cite{thurston}, by another geometric construction, has associated to a triangulation a point $\delta'^\calE$ in the  same lattice  $M^\calE$  (see section \ref{sectthurston}, esp. Thm. \ref{thmthurston} and Cor. \ref{corallcock}). Our final result  (Thm. \ref{mainthm2}) is that (up to isometries of $M^\calE$) the two  points $\delta^\calE$ and $\delta'^\calE$ coincide. This is done by identifying a basis (thought of as a basis of cycles) for $M^\calE$ and evaluating (the cocycles)  $\delta^\calE$ and $\delta'^\calE$ on it respectively. The key point of the identification being that the basis that we consider is constructed purely combinatorially and  $\delta^\calE$ and $\delta'^\calE$ measure the same thing: the combinatorial length of the elements of the basis.

\smallskip

A few words about the organization of the paper. We start by introducing the necessary background on Eisenstein lattices (section \ref{defeisenstein}) and by reviewing the results of Thurston \cite{thurston} on the triangulations of $S^2$ (section \ref{sectthurston}). The only notable aspect here is Cor. \ref{corallcock}, where we identify the Eisenstein lattice and the arithmetic group used by Thurston \cite{thurston} (confirming an earlier guess of Allcock \cite[pg. 294]{allcockaut}). In sections \ref{primcoh} and \ref{sectsmooth}, we discuss the Type III $K3$ surfaces and their cohomology. This is mostly standard material (based on Friedman--Scattone \cite{friedmanscattone}) adapted to our particular needs. In section \ref{sectminus1}, we specialize the discussion to minus-one-forms and obtain more precise results on the structure of the lattices involved in our construction. The structure over the Eisenstein integers is 
 introduced in section \ref{secteisenstein}. At this point, we obtain the first half of our main result: by enriching the structure considered by Friedman--Scattone, we associate to a triangulation of $S^2$ an arithmetic invariant, the orbit of a  point in a  Eisenstein lattice  (Thm. \ref{mainthm1}). Finally, in section \ref{sectequivalence} (Thm. \ref{mainthm2}), we conclude, that under an appropriate identification, our arithmetic invariant is the same as that associated by  Thurston \cite[Thm. 0.1]{thurston}. 

\smallskip

As already noted, it is hoped that this work might lead to a better  understanding of the Type III degenerations of $K3$ surfaces, with the goal of obtaining a  geometric (toroidal) compactification for the moduli space of polarized $K3$ surfaces.
In our opinion there are two  main points made by our paper that make this goal more likely. First, comparing with the results of Friedman--Scattone \cite{friedmanscattone}, it is clear that one gets a richer picture by taking into account polarization for the Type III $K3$ surfaces. Secondly and more importantly,  the results of Thurston can be interpreted as saying that the discrete data of a degeneration of $K3$ surfaces is of arithmetic nature, as needed for a toroidal compactification. Still, the natural question that arises then is if it is possible to  extend the results of Thurston to general Type III $K3$ surfaces (not necessary in minus-one-form) and to interpret the elementary modifications in arithmetic terms.

\subsection*{Acknowledgments} I would like to thank Bob Friedman and Daniel Allcock for some helpful comments early on. While preparing this manuscript, I had discussions on closely related topics with Valery Alexeev, Paul Hacking and Sean Keel. 
 
%%%%%%%%%%%%%%%%%%%%%%%%%%%%%%%%%%%%%%%%%%%%%%%%%%%%%%%%%%%%%%
%%% S1
\section{Eisenstein lattices}\label{defeisenstein} 
By a lattice we understand a free $\bZ$-module together with a symmetric bilinear form.  An Eisenstein lattice is a free module over the Eisenstein integers $\calE=\bZ[\omega]$ together with a hermitian form.  An Eisenstein lattice is called {\it integral} if the corresponding hermitian form is $\calE$ valued (see \cite[pg. 53]{conway}). It is well known that an Eisenstein lattice is equivalent to a $\bZ$-lattice (of double rank)  together with a choice of fixed-point-free isometry of order $3$ (see \cite[\S2.6]{conway}). However, in order to preserve the integrality condition, it is  necessary to use a scaling factor in this equivalence. We explain this below and introduce some notations and results needed later.

\smallskip

Recall first that the data of a complex vector space $W$ together with a hermitian form $H$ is equivalent to  the data of a real vector space together with a symmetric bilinear form $\langle,\rangle$ and a compatible complex structure $J$ (i.e. an isometry with $J^2=-\textrm{id}$). Namely, the restriction of scalars makes $W$ a real vector space and the multiplication by $i$ defines a complex structure $J$.  The bilinear form is given by 
$$\langle x,y\rangle=\mathop{Re} H(x,y).$$
 Conversely, $J$ defines the multiplication by $i$ and the hermitian form is:
$$H(x,y)=\langle x,y\rangle+i\langle x,J(y)\rangle.$$

\smallskip

  Similarly,  a free module $L^\calE$ over the Eisenstein integers is equivalent to a $\bZ$-module $L$ endowed with a fixed-point-free automorphism $\rho$ of order $3$. A hermitian form $H$ on $L^\calE$ determines, as above, a bilinear form $\langle,\rangle$ on $L$, with $\rho$ becoming an isometry. Conversely, $\rho$ defines the multiplication by $\omega=e^{\frac{2\pi i}{3}}$, which makes $L$ an Eisenstein module. Also, $\rho$ defines a complex structure $J$ on $L_\bR=L\otimes_\bZ\bR$ by 
$$J(x)=\frac{1}{\sqrt{3}}\left(\rho(x)-\rho^2(x)\right).$$ 
Thus, we obtain a hermitian form on $L^\calE$: 
$$H(x,y)=\langle x,y\rangle+\frac{i}{\sqrt{3}}\langle x,\rho(y)-\rho^2(y)\rangle.$$
 However, given that $\langle,\rangle$ on $L$ is $\bZ$-valued, it is not necessarily true that $H$ is $\calE$-valued. For this reason,  it is convenient to scale $H$ as follows:
\begin{equation}\label{eqherm}
h(x,y)=\frac{3}{2}H(x,y)=\frac{1}{2}\left(3\langle x,y\rangle+\theta\langle x,\rho(y)-\rho^2(y)\rangle\right),
\end{equation}
where $\theta=i\sqrt{3}=\omega-\omega^2$ (compare with \cite[Eq. (2.2)]{act}). In fact, we even have (see \cite[Lemma 2.1]{act}):
\begin{equation}\label{cond}
h(x,y)\in \theta \calE \textrm{ for all } x,y\in L^\calE
\end{equation}
[N.B. this condition can be viewed as the analogue over $\calE$ of the even lattice condition]. For further reference, we note:
\begin{eqnarray}\label{eqbil}
\langle x,y\rangle&=&\frac{2}{3}\mathop{Re} h(x,y)\\
|x|&=&\frac{2}{3} ||x||\label{eqnorm}
\end{eqnarray}
where the norms are defined by  $|x|=\langle x,x\rangle$ and $||x||=h(x,x)$ respectively.  

\begin{convention}\label{convention}
In what follows, when passing from $\bZ$-lattices to $\calE$-lattices, we always assume that the induced Hermitian form is given by (\ref{eqherm}). With such a scaling, starting with an integral $\bZ$-lattice the associated Eisenstein lattice is integral and satisfies (\ref{cond}). Conversely, the bilinear form induced from a hermitian form on an Eisenstein module is always assumed to be given by (\ref{eqbil}). Assuming in addition that  (\ref{cond}) is satisfied the resulting $\bZ$-lattice is $\bZ$-valued and even.
\end{convention}

Some standard examples of Eisenstein lattices are obtained by putting an Eisenstein  structure on classical $\bZ$-lattices:
 \begin{example}\label{deflat} The root lattices $A_2$, $D_4$, $E_6$, and $E_8$ can be given an Eisenstein lattice structure (see \cite[pg. 3]{artebani3}). This corresponds geometrically to the fact (useful in the study of cubic threefolds \cite{act}) that the singularities $A_2$, $D_4$, $E_6$, and $E_8$ are cyclic triple suspensions of the  singularities $A_1,\dots,A_4$ respectively. Since, up to conjugacy, there exists a unique order $3$ element in the corresponding Weyl group that acts without fixed points (see \cite{carter}), there is no ambiguity in the choice of Eisenstein structure. We denote the resulting Eisenstein lattices  by using the superscript $\calE$ (e.g. $E_8^\calE$). Similarly, the natural analogue of the hyperbolic plane for Eisenstein lattices is $H^\calE:=\left(\begin{array}{cc}0&\theta\\ \bar{\theta} &0\end{array}\right)$, with underlying $\bZ$-lattice $U\oplus U$ (two copies of the hyperbolic plane).
\end{example}

It is well known that there exists a unique  indefinite unimodular  even lattice with a given signature.  A similar fact holds for Eisenstein lattices (\cite[Lemma 2.6]{basak}):
\begin{proposition}
There is at most one indefinite Eisenstein lattice $L$ in a given signature satisfying $L=\theta L^*$ (where $L^*$ is the dual lattice\footnote{The dual lattice is defined by $L^*=\left\{x\in L\otimes_\calE\bQ(\calE)\mid h(x,y)\in\calE \textrm{ for all } y\in L\right\}$. Assuming $L$ is integral, $L\subseteq L^*$. The condition (\ref{cond}) is then equivalent to $\frac{1}{\theta}L\subseteq L^*$ (or $L\subseteq \theta L^*$).}). \qed
 \end{proposition}
 
 An example of a lattice satisfying the assumptions of the previous proposition is the rank $10$ Eisenstein lattice $E_8^\calE\oplus E_8^\calE\oplus H^\calE$  with underlying $\bZ$-lattice $E_8^{\oplus 2}\oplus U^{\oplus 2}$. This lattice was investigated by Allcock \cite{allcockaut}, who related it to the work of Thurston \cite{thurston} (see \cite[pg. 294]{allcockaut}). We arrive to the same Eisenstein lattice by starting with the lattice  $E_8^{\oplus 2}\oplus U^{\oplus 2}$ and endowing it with an Eisenstein structure. The following result says this can be done in a unique way (see also \cite[Cor. 5.3]{artebani3}). 

 \begin{corollary}\label{corunique}
 Let $M=E_8^{\oplus 2}\oplus U^{\oplus 2}$ and assume $\rho\in O(M)$ is a fixed-point-free isometry of order $3$. Then, the associated Eisenstein lattice $M^\calE$ is isometric to $E_8^\calE\oplus E_8^\calE\oplus H^\calE$.
\end{corollary}
\begin{proof} Since $M^\calE$ is constructed from $M$ using the scaling  (\ref{cond}),  we get $M^\calE\subseteq \theta (M^\calE)^*$. For the converse inclusion, we note that
$$\det M^\calE=\pm \sqrt{3^n\cdot \det M}=\pm  \sqrt{3^n} \textrm{ (since $M$ is unimodular)},$$
where $n=\rank_\calE M^\calE=10$ (compare \cite[Eq. (67) on pg. 54]{conway}). Thus,
$$\det (M^\calE)^*=\frac{1}{\det M^\calE}=\frac{\det M^\calE}{3^n}=\det \left(\frac{1}{\theta} M^\calE\right) \textrm{(N.B. $\theta.\bar{\theta}=3$)}.$$
Since $\frac{1}{\theta} M^\calE\subseteq  (M^\calE)^*$, the equality must hold. We conclude, by the previous proposition, that the isometry class $M^\calE$ is uniquely determined. Since $M$ can be obtained from  $E_8^\calE\oplus E_8^\calE\oplus H^\calE$ by forgetting the Eisenstein structure, the conclusion follows.
\end{proof}  

\begin{notation}
In what follows, we use $M$ to denote the lattice $E_8^{\oplus 2}\oplus U^{\oplus 2}$ and $M^\calE$ for the (associated) Eisenstein lattice $E_8^\calE\oplus E_8^\calE\oplus H^\calE$. The groups $\Gamma$ and $\Gamma^\calE$ are the associated isometry groups.
\end{notation}

With our normalization convention, a root in a $\bZ$-lattice $L$ (i.e. an element of norm $-2$) becomes an element of norm $-3$ in an associated Eisenstein lattice $L^\calE$. Thus, it is natural to define:
\begin{definition}
Let $L^\calE$ be an Eisenstein lattice satisfying the condition (\ref{cond}). An element $z\in L^\calE$ of norm $-3$ is called a {\it root}. A root $z$ defines a complex reflection (i.e. an isometry fixing pointwise a hyperplane) of order $3$ by:
\begin{equation}
s_z(x)=x-(1-\omega)\frac{h(x,z)}{h(z,z)}\cdot z
\end{equation} 
Note that $s_z$ sends $z$ to $\omega\cdot z$ and fixes pointwise the hyperplane $z^\perp$. When viewed over $\bZ$, a root in an Eisenstein lattice corresponds to the root system $A_2=\Span_\bZ(z,\omega\cdot z)$. 
\end{definition}

We close by noting a few arithmetic facts about the lattice $M^\calE:=E_8^\calE\oplus E_8^\calE\oplus H^\calE$ mostly due to Allcock \cite{allcockaut}. First, it is easy to see that $M^\calE$ is generated by $10$ roots forming an $A_{10}$ type graph. Specifically, there exist $\{z_1,\dots,z_{10}\}$ forming a $\calE$-basis for $M^\calE$ such that: 
\begin{equation}\label{herm1}
h(z_k,z_l)=\left\{
\begin{array}{ll}
0 &\textrm{if } k<l-1\\
\bar{\theta} &\textrm{if } k=l-1\\
-3 &\textrm{if } l=k\\
\end{array}
\right. .
\end{equation}
\begin{remark}
It is interesting to note that the $\bZ$-lattices underlying the Eisenstein lattice generated by roots forming an $A_1$, $A_2$, $A_3$ or $A_4$ graph are $A_2$, $D_4$, $E_6$, or $E_8$ respectively (compare with \ref{deflat}). The $A_5$ graph corresponds to negative semi-definite lattice $\widetilde{E}_8$ 
\end{remark}

Furthermore, Allcock \cite[Thm. 5.1]{allcockaut} shows that the group of isometries of $M^\calE$ is generated by reflections in roots. In fact, $11$ complex reflections suffice, corresponding to the fact that $\Gamma^\calE$ can be related to a representation of the spherical braid group on $12$ strands, with the generators mapping to order $3$ reflections. In fact, Allcock \cite[\S5]{allcockaut} notes that the projectivized group $\bP\Gamma^\calE\subset \bP U(1,9)$ coincides with the Deligne--Mostow group $\Gamma_\mu$ with $\mu=(\frac{1}{6},\dots,\frac{1}{6})$. Also, cf. \cite[pg. 294]{allcockaut},  $\Gamma^\calE$ is the discrete group occurring in Thurston \cite{thurston}. Finally, let $\Gamma$ denote the subgroup of isometries of the lattice $E_8^{\oplus 2}\oplus U^{\oplus 2}$ that preserve the two components of the associated type IV domain, and $\rho\in \Gamma$ an order $3$ fixed-point-free isometry. Then, the associated Eisenstein lattice is isometric to $M^\calE\cong E_8^\calE\oplus E_8^\calE\oplus H^\calE$ (cf. \ref{corunique}) and the two groups $\Gamma^\calE\subset U(1,9)$ and   $\Gamma\subset O(2,18)$ are related by:
\begin{equation}\label{mongroup}
\Gamma^\calE=\{\gamma\in \Gamma\mid \gamma\cdot\rho=\rho\cdot \gamma\}
\end{equation}

%%%%%%%%%%%%%%%%%%%%%%%%%%%%%%%%%%%%%%%%%%%%%%%%%%%%%%%%%%%%%%
%%% S2
\section{Summary of the results of Thurston on triangulations of $S^2$}\label{sectthurston}
The purpose of this paper is to discuss the following result of W. Thurston  \cite[Thm. 1]{thurston} regarding the triangulations of $S^2$ from the perspective of degenerations of $K3$ surfaces. Here we briefly review the construction of Thurston and see that the $K3$ surfaces occur naturally in this picture. 

\begin{theorem}[Thurston \cite{thurston}]\label{thmthurston} There exists an Eisenstein lattice $N^\calE$ of signature $(1,9)$ and a group of lattice isomorphisms $\Gamma^\calE$ such that: to every triangulation $T$ of $S^2$ of non-negative combinatorial curvature one can associate a point $\delta^\calE\in N^\calE$ of positive norm, equal to $t$, the number of triangles  in $T$, in such a way that the orbit of $\delta^\calE$ (w.r.t. $\Gamma^\calE$) is a complete invariant of $T$. 
\end{theorem}

We recall that by a triangulation of non-negative combinatorial curvature, one understands a triangulation such that the degree of each vertex of the triangulation is at most $6$. We note that the notion of  triangulation used here allows the identification of two edges of the same triangle (see \cite[pg. 517]{thurston}). 

\smallskip

The lattice $N^{\calE}$ of the theorem is not explicitly identified in \cite{thurston}, but Allcock  \cite[pg. 294]{allcockaut} suggests that, after scaling by $\frac{3}{2}$, $N^{\calE}$ is isometric to  $M^\calE=E_8^\calE\oplus E_8^\calE\oplus H^\calE$.  This is based on the fact, already mentioned above,  that the discrete group occurring in Thm. \ref{thmthurston} coincides with the isometry group $\Gamma^\calE$ of $M^\calE$ (thus, we use the same notation).  As explained in  \S\ref{sectdm} below,  indeed $N^\calE=M^\calE(\frac{2}{3})$ (see Cor. \ref{corallcock}).   Also, we note that  implicitly contained in Thurston \cite{thurston} is the fact that the sublattice $R^\calE$ spanned by the roots contained in the orthogonal complement $\langle\delta^\calE\rangle^\perp$ is a direct sum of $A_2^\calE$, $D_4^\calE$, $E_6^\calE$, and $E_8^\calE$ summands, that determines the degree type of the triangulation (see \S\ref{sectdegtype}). 

\begin{remark}\label{convention2} In order to preserve the integrality condition (see \S\ref{defeisenstein}), we prefer to work with the scaled version $M^\calE= N^\calE\left(\frac{3}{2}\right)$. Thus, in what follows we use $M^\calE$ (and $\delta^\calE\in M^\calE$, etc.) instead of $N^\calE$ and assume the scaling convention \ref{convention} in place. It is important to note that, using the scaling or not, the $\bZ$-lattice corresponding to the rank $1$  sublattice generated by the element $\delta^\calE\in N^\calE$ of Thm. \ref{thmthurston} is $A_2\left(-\frac{t}{2}\right)$. 
 \end{remark}

\subsection{The construction of Thurston}\label{sectdm}
 The approach of Thurston to Thm. \ref{thmthurston} is to associate to a triangulation a geometric object and then to consider the moduli of those objects. Specifically, by declaring every triangle in the triangulation $T$ an equilateral triangle, one obtains a flat Euclidean metric on $S^2$ with cone type singularities. By Euler formula, for a generic triangulation of non-negative curvature, there are $12$ singular points of prescribed cone curvatures. Then, up to scaling the metric, the geometric data associated to $T$ is equivalent to specifying $12$ points $\{p_1,\dots,p_{12}\}$ in $\bP^1$. By work of Deligne--Mostow, the moduli of $12$ points in $\bP^1$ can be uniformized by a ball. Thus, from a algebro-geometric point of view the construction of \cite{thurston}, can be viewed as the following association:
$$
\begin{array}{lclclclcl}
&\textrm{\cite[\S0]{thurston}}&&\textrm{\cite[\S8]{thurston}}&& \textrm{\cite[\S3]{thurston}}\\
T& \Longrightarrow& (S^2,g)&  \Longrightarrow& (\bP^1,(p_1,\dots,p_{12}))
&\Longrightarrow& [\delta^\calE]\in \calB/\bP\Gamma^\calE\\
 \textrm{triang.}& &\textrm{cone metric}&& \textrm{pointed } \bP^1&&\textrm{periods}\\
\end{array}$$
where $\calB$ is the $9$-dimensional complex ball:
\begin{equation}\label{ball}
\calB=\left\{[z]\in\bP(N^\calE\otimes_{\calE}\bC)\mid ||z||>0\right\},
\end{equation}
where $[z]=z\cdot \bC$ and $||\cdot||$ is norm w.r.t. to the hermitian form induced from $N^\calE$. The  periods corresponding to triangulations belong to the discrete set $\bP(N^\calE)\cap \calB_9$. Due to the scaling of the metric, the knowledge of the period point $[\delta^\calE]$ (defined only up to scaling by $\bC^*$) only determines the shape of the triangulation. Thurston shows that in fact one can choose a lift $\delta^\calE\in N^\calE$ in such a way that one recovers the triangulation. 

\begin{remark}\label{remreco}
A few words about the reconstruction procedure. First, by a Torelli type theorem, the period point determines the cone metric up to scaling. Choosing a norm $t$ representative $\delta^\calE\in N^\calE\otimes_{\calE}\bC$ for the period, gives the correct scaling of the metric. Then,  the cone metric  determines a Delaunay triangulation of $S^2$ with vertices at the cone points (see  Thurston \cite[Prop. 3.1]{thurston} and Rivin \cite[\S10]{rivin}). Using this, one opens the sphere to get a star-like polygon in the Euclidean plane (see \cite[\S7]{thurston}, esp. \cite[Prop. 7.1]{thurston} and  \cite[Fig. 14]{thurston}).
To recover the triangulation, one needs  a reduction of the structure to the Eisenstein integers
(see \cite[pg. 535]{thurston}), i.e. an isometric embedding of the Eisenstein lattice into the plane such that the vertices of the polygon are points in the lattice. By using a vertex of the polygon as the origin, we can identify the Euclidean plane with  $\bC$ and the Eisenstein lattice with $\calE\subset \bC$, but there is still ambiguity due to rotations, given by $S^1\cap \bQ(i\sqrt{3})$, that place the second vertex into the lattice $\calE$. In terms of periods, it means that we need a representative $\delta^\calE\in N^\calE$ for the period. A choice $\delta^\calE\in N^\calE\otimes_\calE \bQ(i\sqrt{3})$, say of correct norm, is typically not enough to recover the triangulation.
 \end{remark}

There are two points that need to be explained in the above identification. First, we briefly explain the relation to the work of Deligne--Mostow \cite{delignemostow}. It is clear that the space of cocycle of Thurston (see \cite[\S3]{thurston}) is the cohomology group $\sH^1(U,\bL)$, where $U=\bP^1\setminus\{p_1,\dots,p_{12}\}$ and $\bL$ is the $1$-dimensional local system associated to the orthogonal holonomy representation $H_0:\pi_1(U,x)\to S^1\subset \bC^*$ given by mapping an oriented loop  around a singular point to $e^{\frac{5\pi i}{3}}$ (see \cite[pg. 526]{thurston}). The cocycle space  $\sH^1(U,\bL)$ is a basic ingredient in Deligne--Mostow \cite{delignemostow} and, following the exposition of Looijenga \cite{loijl} (see esp. \cite[Rem. 1.1]{loijl}), it is  not hard to see that  indeed $[\delta^\calE]=\delta^\calE\cdot \bC$ is the period point for the corresponding $12$ points in $\bP^1$. Secondly, we need to explain the Eisenstein structure of the space of cocycles. It is known that $\sH^1(U,\bL)$ can be naturally defined over the ring of Eisenstein integers (see \cite[\S12]{delignemostow} and \cite[\S4.3]{loijl}). The cocycles with $\calE$ coefficients are precisely those coming from triangulations (cf. \cite[pg. 535]{thurston}) [N.B. in particular, the periods for points coming from triangulations are those defined over the imaginary quadratic field $\bQ(i\sqrt{3}$)].

\subsection{The Eisenstein lattice $N^\calE$ of Thm. \ref{thmthurston}}\label{identn} 
To understand the Eisenstein lattice structure on the space of cocycles with Eisenstein coefficients, we note that an algebro-geometric interpretation of this space can obtained by considering the cyclic order $6$ cover $C$ of $\bP^1$ totally branched at the $12$ special points, i.e. 
$$C: z^6=f_{12}(t),$$
where $f_{12}(t)$ is a degree $12$ polynomial defining the $12$ points $\bP^1$ (see \cite[(2.23)]{delignemostow} and \cite[\S4]{loijl}). Then, there exists a natural identification: 
\begin{equation}
\sH^1(U,\bL)=\sH^1(C,\bC)_{-\omega},
\end{equation}
 where the eigenspace is taken with respect to the $\mu_6$-action induced on cohomology by the automorphism $\psi:C\to C$ given by $\psi(z,t)=(\zeta\cdot z,t)$, where $\zeta=-\overline{\omega}=e^{\frac{2\pi i}{6}}$ is an order $6$ primitive root of unity.
From this perspective, the period point associated to $12$ points is simply the line: 
$$\sH^{1,0}(C,\bC)_{-\omega}\subset \sH^1(C,\bC)_{-\omega}.$$
Note that $\sH^1(C,\bC)_{-\omega}$ is a $10$-dimensional complex vector space, that carries a natural Hermitian form of hyperbolic signature $(1,9)$. The subspace  $\sH^{1,0}(C,\bC)_{-\omega}$ is $1$-dimensional and positive definite (see \cite[Ex. 3.11]{vangeemenhalf} for an explicit proof). Everything can be defined over the Eisenstein integers by replacing $\sH^1(C,\bC)$ by $\sH^1(C,\calE)$. We conclude that the Eisenstein lattice $N^\calE$ mentioned in Thm. \ref{thmthurston} is $\sH^1(C,\calE)_{-\omega}\cong (\calE)^{10}$ endowed with a scaling of the natural Hermitian form induced from the cup product on $\sH^1(C,\calE)$. A computation of the Hermitian form was done by Looijenga \cite[Prop. 4.6]{loijl}. Putting everything together, we obtain:
\begin{corollary}\label{corallcock}
The Eisenstein lattice is $N^\calE$ occurring in Thm. \ref{thmthurston} is a scaling by $\frac{2}{3}$ of the Eisenstein lattice $M^\calE=E_8^\calE\oplus E_8^\calE\oplus H^\calE$. Also, the group $\Gamma^\calE$ of the Thm. \ref{thmthurston} is the isometry group of $M^\calE$ and is generated by complex reflections. In particular,  the underlying $\bZ$-lattice (obtained by taking the real part of the hermitian form) of  $N^\calE$  is $E_8^{\oplus 2}\oplus U^{\oplus 2}$.
\end{corollary}
\begin{proof}
As mentioned before, the statement about $\Gamma^\calE$ is due to Allcock \cite[\S5]{allcockaut}, who also conjectured the structure of $N^\calE$. Here, we explicitly check this. Let $h$ be the Hermitian form on $M^\calE$ and $H$ the Hermitian form considered by Looijenga \cite{loijl}.  Looijenga \cite[Prop. 4.6]{loijl} has computed that the hermitian form $H$ is given by: 
\begin{equation}\label{herm2}
H(\epsilon_k,\epsilon_l)=\left\{
\begin{array}{ll}
0 &\textrm{if } l<k-1\\
-\frac{1}{2} &\textrm{if } l=k-1\\
\frac{\sqrt{3}}{2} &\textrm{if } l=k\\
\end{array}
\right.
\end{equation}
with respect to a certain basis $\{\epsilon_1,\dots,\epsilon_{10}\}$ of $\sH_1(C,\calE)_{-\overline{\omega}}$ (the dual of the space of cocycles). Under an appropriate identification, it is easy to see that the hermitian forms $h$ (see formulas (\ref{herm1})) and $H$ are related by: 
\begin{equation}\label{eqcb1}
H(z,w)=-\frac{1}{2\sqrt{3}}\cdot h(z,w).
\end{equation}
Specifically, one obtains  (\ref{herm1}) from (\ref{herm2}) using the change of basis given by: 
\begin{equation*}
z_k=\overline{i}^k \cdot \epsilon_k
\end{equation*}
Note that although $i\not\in \calE$, the change of basis is well defined over $\calE$. Namely,  $\{\epsilon_1,\dots,\epsilon_{10}\}$ is obtained from an $\calE$-basis $\{\epsilon_1',\dots,\epsilon_{10}'\}$\footnote{with $\epsilon_k'=\overline{w}_k \sum_{g\in \mu_m} \chi(g) g_*\delta_k$ in the notation of  \cite[\S4.3]{loijl}} for $\sH_1(C,\calE)_{-\overline{\omega}}$    by the transformation: 
$\epsilon_k=\overline{w}_k \cdot \epsilon_k'$,
where $w_k=e^{i\pi (\mu_0+\dots+\mu_{k-1})}=e^{\frac{k\pi i}{6}}$. Since $z_k=\overline{i}^k \cdot \overline{w}_k \cdot \epsilon_k'$ and $\overline{i}^k \cdot \overline{w}_k=e^{-\frac{2k\pi i}{3}}\in \calE$, we conclude that the identification (\ref{eqcb1}) is indeed valid over $\calE$. 

The claim about $N^\calE$ now follows from (\ref{eqcb1}) and the fact that the Hermitian form used by Thurston \cite{thurston} on the space of cocycles is $-\frac{4}{\sqrt{3}}H$. The minus sign is due to a sign convention used in \cite{loijl} and $\frac{\sqrt{3}}{4}$ is the area of a standard equilateral triangle. Thus, in Thurston \cite{thurston}, to pass from the total area to the number of equilateral triangles, one has to scale the Hermitian form by $\frac{4}{\sqrt{3}}$. 
\end{proof}

\begin{remark}\label{remk3} 
We note that one can associate naturally to (a shape of) a triangulation of $S^2$ a smooth $K3$ surface $X$. Namely, to a triangulation of $S^2$ one associates $12$ points in $\bP^1$ and then a  $K3$ surfaces by considering the elliptic fibration $X\to \bP^1$:
\begin{equation}\label{eqelliptic}
X/\bP^1:\ y^2=x^3+f_{12}(t).
\end{equation}
This procedure is known to give an equivalent construction of the Deligne--Mostow uniformization of the moduli space space of $12$ points in $\bP^1$ (see  van Geemen \cite[Ex. 3.11]{vangeemenhalf} and Kondo \cite[\S4]{k2}). A key ingredient  being that  $X$ has a non-symplectic automorphism $\sigma$ of order $3$ (i.e. $\sigma(x,y,t)=(\omega\cdot x,y,t)$). One can then consider the eigenperiods of $X$ (see Dolgachev--Kondo \cite[\S11]{dk}) and obtain a ball quotients uniformization for the moduli of $12$ points in $\bP^1$.
\end{remark}

%%%%%%%%%%%%%%%%%%%%%%%%%%%%%%%%%%%%%%%%%%%%%%%%%%%%%%%%%%%%%%
%%% S3
\section{The lattice associated to a Type III $K3$ surface}\label{primcoh}
In this section, we associate to a type III $K3$ surface $X_0$ an even lattice $\overline{L}$ of signature $(1,18)$, that represents the interesting cohomology of $X_0$. In the context of our paper, the lattice $\overline{L}$ should be viewed as a discrete arithmetic invariant associated to a Type III $K3$ surface. Although at various points we make use of analytic results, we emphasize that $\overline{L}$ is constructed purely from the  combinatorial data of $X_0$. Most of the material in this and following section is based on Friedman--Scattone \cite{friedmanscattone}. Also related is M. Olsson \cite[\S3]{olsson}.
 
\smallskip 

To fix the notation, we recall that a (combinatorial) Type III $K3$ surface $X_0$ is  a normal crossing variety
$$X_0=\cup_{i=1,\dots,n} V_i$$ 
 such that the components $V_i$ are rational surfaces glued along anti-canonical cycles and the dual graph of $X_0$ is a triangulation $T$ of $S^2$. We denote $\widetilde{V}_i$ the normalizations of the components,  $D_{ij}$  the double curves, and $t_{ijk}$ the triple points. The double curves are assumed to satisfy the triple point formula:
\begin{equation}\label{triplepoint}
D_{ij}^2+D_{ji}^2=-2
\end{equation}
(with the convention that $D_{ij}$ and $D_{ji}$ lie on $\widetilde{V}_i$ and $\widetilde{V}_j$ respectively), except the nodal case, when say $D_{ij}$ is nodal on $\widetilde{V}_i$, in which case we require $D_{ij}^2+D_{ji}^2=0$. 

\begin{notation}\label{notations1}
We use $n$, $e$, and $t$ to denote the number of components, of double curves, and of triple points of $X_0$ respectively. For the dual graph $T$, $n$ is the number of vertices and $t$ is the number of triangles. In particular, from the Euler relation $n-e+t=2$, it follows that $e=3n-6$ and $t=2n-4$. Also, $C_i$ and $C^i$ denote the chains and cochains of the triangulation $T$ (e.g. $C^1$ is an abelian group isomorphic to  $\bZ^e$ and with natural generators corresponding to the double curves $D_{ij}$; its dual, $C_1$ has generators corresponding to the edges of the triangulation). For consistence, we mention that an orientation of $T\cong S^2$ is fixed throughout and the $1$-simplices in $T$ are oriented from the vertex $v_i$ (corresponding to $V_i$) to the vertex $v_j$ for $i>j$.
\end{notation}

\subsection{The cohomology of a Type III $K3$ surface}\label{sectcoho}  
The cohomology of the normal crossing variety $X_0$ is computed by the Mayer-Vietoris spectral sequence:
$$E_1^{pq}=\sH^q(X_0^{[p]})\Longrightarrow \sH^{p+q}(X_0)$$
In particular, we obtain the exact sequence:
\begin{equation}\label{computeh2}
0\to (W_0\sH^2)_\bZ\to \sH^2(X_0,\bZ)\to (\Gr_2^W\sH^2)_\bZ\to 0,
\end{equation}
%(see \cite[pg. 16]{friedmanscattone}), 
where $(W_0\sH^2)_\bZ\cong \sH^2(S^2,\bZ)\cong \bZ$ and
\begin{equation}\label{defl}
L:=(\Gr_2^W\sH^2)_\bZ=\ker\left(\oplus_{i=1}^n \sH^2(\widetilde{V}_i,\bZ)\to\oplus_{ij}\sH^2(D_{ij},\bZ)\right)
\end{equation}
Note that the surfaces $\widetilde{V}_i$ are rational surfaces, whose Betti numbers can be computed in terms of the self-intersections of the double curves (N.B. $b_2(V_i)=10-K_{V_i}^2$). By combining the triple point formula and Euler formula,  we get $\rank \oplus_{i=1}^n \sH^2(\widetilde{V}_i,\bZ)=4n+12$. 
The map $\oplus_{i=1}^n\sH^2(\widetilde{V}_i)\to\oplus_{ij}\sH^2(D_{ij})$ is surjective over $\bQ$ (cf. \cite[Prop. 7.2]{friedmanscattone}). 
%(fact related to the non-existence of infinitesimal automorphisms for $X_0$). 
We conclude that rank of $L$ is $18+n$.

\smallskip

The cohomology $\sH^2(X_0)$ carries a natural mixed Hodge structure, which is an extension of a trivial Hodge structures of weight $2$ by a Hodge structure of weight $0$. By Carlson theory, such extensions are classified by a certain extension homomorphism:
$$\psi:L\to \bC^*,$$
where $L=(\Gr_2^W\sH^2)_\bZ$ and $\bC^*\cong(W_0\sH^2)_\bZ\otimes \bC/(W_0\sH^2)_\bZ$ (see \cite[Sect. 3 and 4]{friedmanscattone}). Relevant for us is the fact that the group of Cartier divisors on $X_0$ can be identified with $\ker \psi\subseteq L$. In other words, $L$ represents the cohomology classes that satisfy the combinatorial conditions of Cartier divisor, while those in $\ker \psi$ also satisfy the analytic gluing conditions. In particular, we note that for each component $V_i$ of $X_0$ there is a naturally associated combinatorial Cartier divisor:
$$\xi_i:=\sum_j \left(D_{ij}-D_{ji}\right)\in L.$$
Note the obvious linear relation
\begin{equation}
\sum_1^n \xi_i=0
\end{equation}
satisfied by the $\xi_i$.

\smallskip

By results of Friedman \cite{friedmansmooth}, it is known that $X_0$ is smoothable if and only if $\xi_i$ is indeed a Cartier divisor for each $i$, i.e. $\xi_i\in \ker \psi$ for all $i$ (in which case we say $X_0$ is {\it d-semistable}). Thus, it is natural to consider 
\begin{equation}\label{defk}
K=\Span(\{\xi_i\}_i)\subset L
\end{equation}
 and define 
\begin{equation}\label{defbarl}
\overline{L}:=L/K \textrm{ modulo torsion}.
\end{equation}
From the point of view of degenerations, $\overline{L}$ is the interesting cohomology of $X_0$. For example, assuming a maximal algebraic semistable degeneration $\calX\to\Delta$ of $K3$ surfaces, i.e. $L=\ker \psi$ and the nearby smooth fibers $X_t$ have Picard rank $19$, we get $\Pic(\calX)=\Pic(X_0)=L$. Then, $\overline{L}$ is the group of line bundles on the total space $\calX$ modulo twisting by $\calO_\calX\left(\sum \alpha_i V_i\right)$ .

\begin{remark} We always implicitly assume that $X_0$ is smoothable (or $d$-semistable, i.e.  $K\subseteq \ker \psi$). However, for our purposes, the study of the discrete invariants of a Type III $K3$ surface, one can  ignore this condition: it   
can be always achieved by a locally trivial deformation. On the other hand,  the combinatorics of the surface  is preserved under such deformation. 
\end{remark} 

\subsection{The lattice structure of $\overline{L}$}\label{sectlatl} 
A priori, $\overline{L}$ is only a free abelian group. Here, we show that $\overline{L}$ caries a natural bilinear form of signature $(1,18)$. To start, we define $H$ to be the direct sum lattice:  
$$H:=\oplus\sH_2(\widetilde{V_i},\bZ),$$
where each summand is equipped with the natural intersection form. By Poincare duality on each of the components, it  follows easily that $H$ is an odd, unimodular, lattice of signature $(n,3n+12)$, thus isometric to $(1)^{n}\oplus (-1)^{3n+12}$. In particular, $H$ can identified with its dual $H^*=\oplus\sH^2(\widetilde{V_i},\bZ)$. 

\begin{remark}\label{remint}
This choice for the intersection form on the singular fiber is discussed in \cite[\S2.4]{persson}. Its main property is that it behaves well with respect to smoothings. Namely, if $\calD$ and $\calD'$ are divisors (in a topological sense) on the total space of a semi-stable degeneration then $D_t.D_t'=D_0.D_0'$ (see  \cite[Prop. 2.4.1]{persson}), where $D_t$ and $D_0$ denote the restrictions of $\calD$ (and similarly for $\calD'$) to a smooth fiber  $X_t$ and to the central fiber $X_0$ respectively. 
\end{remark}

We define then two  sublattices $D$  and $L$ of $H$ by:
$$D:=\im\left(\oplus \sH_2(D_{ij},\bZ)\to\oplus \sH_2(\widetilde{V_i},\bZ)\right)\hookrightarrow H,$$
and, as before,  
$$L:=\ker\left(\oplus \sH^2(\widetilde{V}_i,\bZ)\to\oplus\sH^2(D_{ij},\bZ)\right)\hookrightarrow H^*\cong H.$$
Clearly, $L$ can be defined equivalently as the orthogonal complement of $D$ in $H$, i.e. $L=D^{\perp}_H$. Since $D$ and $L$ are  mutually orthogonal sublattices in $H$, they intersect in an isotropic sublattice, which at least over $\bQ$ is the common radical of both. It is clear that 
\begin{equation}\label{equalityk}
K\subseteq L\cap D
\end{equation}
(see the definition (\ref{defk}) of $K$). Since $L$ is primitive by definition (it is the kernel of a map of free abelian groups), the equality in (\ref{equalityk}) holds iff  $K$ is the radical of $D$. This is indeed the case as shown by the following lemma.

\begin{lemma}
With notations as above, $K$ is a primitive isotropic sublattice of $D$ with the property that $D/K$ is a negative definite lattice. In fact, $D/K$ is isometric to the root lattice $A_{t-1}$. In particular, $K$ is the radical of $D$ and thus
$$D\cap L=K.$$
\end{lemma}
\begin{proof}
Let $\overline{D}=D/K$. The statement of the lemma is that the $\bZ$-module $\overline{D}$ is torsion free and the induced intersection form on $\overline{D}$ (well defined for $K$ totally isotropic) is negative definite. We prove these facts by relating the inclusion $K\subset D$ to the combinatorics of the  dual graph $T$ of $X_0$. First, there is a natural identification of $D$ with $C^1$ (see \ref{notations1} for notations and conventions) by means of the map:
$$C^1= \oplus \sH^0(D_{ij})\xrightarrow{P.D.}\oplus \sH_2(D_{ij})\twoheadrightarrow D\subset H=\oplus \sH_2(\widetilde{V}_i)$$
[N.B. writing $(D_{ij}-D_{ji})\in D$ gives an orientation of the edge from $i$ to $j$ in $C_1$, say from $i$ to $j$; $(D_{ji}-D_{ij})$ corresponds to the reverse orientation]. On the other hand, from the perspective of degenerations of $K3$ surfaces, $K$ should be thought  as coming from a mapping of $\sH_4(X_0)=\oplus \sH_4(\widetilde{V_i})\cong \sH^0(\widetilde{V_i})=C^0$ into $\sH^2(X_0)$. Purely combinatorially, this means that $K\subset D$ is the image of the composite map $C^0\to C^1\cong D$. This  is easily checked by noting that the natural generators of $C^0$ map to the generators $\xi_i$ of $K$. Thus, we obtain the  following commutative diagram of $\bZ$-modules:
\begin{equation}\label{eqdbar}
\xymatrix{
0\ar@{->}[r]&\sH^0(S^2)\ar@{->}[r]&C^0\ar@{->>}[d]\ar@{->}[r]&C^1\ar@{->}[r]\ar@{..>>}[rd]\ar@{=}[d]&C^2\ar@{->}[r]&\sH^2(S^2)\ar@{->}[r]&0\\
&0\ar@{->}[r]&K\ar@{->}[r]\ar@{^{(}..>}[ru]&D\ar@{->}[r]&\bar{D}\ar@{^{(}->}[u]\ar@{->}[r]&0
}
\end{equation}
Since $\overline{D}$ can be identified with the image of $C^1\to C^2$, $\overline{D}$ is torsion free. 

To prove that $\overline{D}$ is  negative definite, we note that the identification $\overline{D}=\im(C^1\to C^2)$ becomes a lattice isometry when we endow $\im(C^1\to C^2)$ with the obvious lattice structure. Namely, identify $C^2$ with $\bZ^t$ with the lattice structure $(-1)^t$ by declaring the generators of $C^2$ to be of norm $-1$. Then, for appropriate choice of  orientation,   $\im(C^1\to C^2)$ is identified to $\left\{(x_1,\dots,x_t)\in \bZ^t\mid x_1+\dots+x_t=0\right\}$, the standard $A_{t-1}$ lattice. Geometrically, the roots of $A_{t-1}$ can viewed are differences $t_\alpha-t_\beta$, where $t_\alpha$ and $t_\beta$ correspond to the triple points of $X_0$. At this point, we only have a linear isomorphism $\phi$ between two lattices $\overline{D}$ and $ A_{t-1}\cong\im(C^1\to C^2)$. To get an isometry, we need to check that for some set of generators  of $\overline{D}$ we get the expected intersection numbers. We consider the obvious set of generators for $\overline{D}$: the classes $[D_{ij}-D_{ji}]$ of $(D_{ij}-D_{ji})\in D$. The identification $\phi$ regards $[D_{ij}-D_{ji}]$ as the difference of the two triple points on the double curve $D_{ij}$. Thus, $\phi( [D_{ij}-D_{ji}])$ is a root of $A_{t-1}$. Indeed, we get the expected intersection number:
$$[D_{ij}-D_{ji}]^2=(D_{ij}-D_{ji}).(D_{ij}-D_{ji})=D_{ij}^2+D_{ji}^2=-2$$
by the triple point formula. Note that in the exceptional nodal case, when $D_{ij}^2+D_{ji}^2=0$, the class $[D_{ij}-D_{ji}]$ is the class of $\xi_i$ or $-\xi_j$, and thus $0$ in $\overline{D}$. Similarly, the intersection of two different classes $[D_{ij}-D_{ji}]$ and $[D_{kl}-D_{lk}]$ is $0$, unless $\{i,j,k,l\}$ is the index set of a triple point, i.e. either $i$ or $j$ coincides with $k$ or $l$. Say $i=k$, we have 
\begin{equation}\label{eqpm}
[D_{ij}-D_{ji}].[D_{il}-D_{li}]=D_{ij}.D_{il}=1.
\end{equation}
This corresponds to the standard fact that two roots $t_\alpha-t_\beta$ and $t_\gamma-t_\delta$ in $A_{t-1}$ intersect  if and only if they share one of the triple points $t_*$. If two distinct (up to $\pm$) roots intersect, the intersection number is either $1$ or $-1$. Via the identification $\phi$, it is easily seen that the intersection number in (\ref{eqpm}) has the right sign. We conclude, as needed, that $\overline{D}$ is isometric to $A_{t-1}$.
\end{proof}

\begin{remark}
As noted by the proposition, $D$ and $\overline{D}$ have very simple lattice structures. Essentially, the only combinatorial information contained in them is the number of triple points. As we will see later, the same is true for $L$.
\end{remark}

Since $K$ is a totally isotropic subspace of $L$, the quotient module
$$\overline{L}=L/\Sat(K)$$
carries a natural intersection form induced from  the lattice $H$. We recall that 
$$\Sat(K)=\left\{x\in L\mid n\cdot x\in K \textrm{ for some } n\in \bZ\setminus\{0\}\right\}$$
and note that the definition of $\overline{L}$ coincides with the earlier one (\ref{defbarl}).

\begin{definition}\label{defpol}
Let $X_0$ be a combinatorial Type III $K3$ surface. By the (cohomology) {\it lattice associated} to $X_0$ we understand the lattice $\overline{L}$ defined as above. By a {\it polarization} for $X_0$ we understand an element $h\in \overline{L}$ with $h^2>0$ (w.r.t. the intersection form on $\overline{L}$).
\end{definition} 

\begin{remark}
The lattice $\overline{L}$ and the notion of polarization for a Type III $K3$ surface appear previously in literature, see for example Friedman--Scattone \cite{friedmanscattone} and M. Olsson \cite[Def. 1.2 and \S3.9, esp. Prop. 3.10]{olsson}. %In fact, $\overline{L}=\sH^1(X,\calM_X^{gp})$ in the notation of \cite{olsson} (see esp. \cite[(3.2.3)]{olsson}).
\end{remark}

By construction, the intersection form on $\overline{L}$ is non-degenerate (N.B. $\Sat(K)$ is the radical of $L$). We conclude  by noting that Hodge index theorem holds for $\overline{L}$:
\begin{proposition}\label{lemmaeven}
The lattice $\overline{L}$ associated to a Type III $K3$ surface is an even lattice of signature $(1,18)$. 
\end{proposition}
\begin{proof} Since $\rank L=18+n$, it follows that $\overline{L}$ has rank  $19$. Let $(a,19-a)$ be the signature of $\overline{L}$. Since  $H$ has signature $(n,3n+12)$ and $K^\perp_H=\overline{L}\oplus \overline{D}\oplus K$ (over $\bQ$) has codimension $n-1$, we get $a\ge 1$. Conversely, let $h_i\in \sH^2(\widetilde{V_i},\bZ)$ be the class of an ample divisor for the component $V_i$ of $X_0$ (for $i=1,\dots,n$). The linear space spanned by $n-1$ of these classes does not meet $L$. Namely, say $\ell=\sum_{i=1}^{n-1} \alpha_i h_i$ belongs to $L$, i.e. $\ell.D_{ij}=\ell.D_{ji}$ for all $i, j$. The ampleness of $h_i$ implies that $h_i.D_{ij}\neq 0$ for all $j$.  Thus, if $h_i$ occurs with non-zero coefficient in $\ell$, the same must be true for all $h_j$ with $V_j$ adjacent to $V_i$. By assumption, $h_n$ does not occur in $\ell$, a contradiction.  Since $\Span(h_1,\dots,h_{n-1})\cap L=0$ and  $\Span(h_1,\dots,h_{n-1})$ is positive definite, we conclude that $\overline{L}$ has signature $(1,18)$.

Let $\ell \in L$, and write $\ell=\sum_i \ell_i$ with $\ell_i\in \sH^2(\widetilde{V_i},\bZ)$. By adjunction formula,  $\ell_i.(\ell_i+K_{\widetilde{V_i}})$ is even. Since $\ell\in L$ and $K_{\widetilde{V_i}}=-\sum_i D_{ij}$, we also get
$$\ell.\sum_i K_{\widetilde{V_i}}=-\ell.\sum_{i,j}  D_{ij}=-2\ell.\sum_{i<j}D_{ij}\equiv 0 \mod 2$$ 
and then
$$\ell.\ell\equiv \ell.\left(\ell+\sum K_{\widetilde{V_i}}\right)\equiv \sum_i \ell_i.(\ell_i+K_{\widetilde{V_i}})\equiv 0 \mod 2.$$
Thus, $\overline{L}$ is an even lattice.  
\end{proof}

%%%%%%%%%%%%%%%%%%%%%%%%%%%%%%%%%%%%%%%%%%%%%%%%%%%%%%%%%%%%%%
%%% S3
\section{Smoothings of Type III degenerations}\label{sectsmooth}
If a type III $K3$ surface $X_0$ is d-semistable, it can be viewed as the central fiber of a semi-stable family $\calX\to \Delta$. By investigating the behavior of the cohomology with respect to the degeneration, we give a precise relation between the lattice  $\overline{L}$ associated to $X_0$ (Def. \ref{defpol}) and the $K3$ lattice $\Lambda:=E_8^2\oplus U^{\oplus 3}$. As before, we are interested only in the discrete aspects of the degeneration, thus we are working with $\bZ$ coefficients and essentially ignore the analytic aspects (Hodge filtration, etc.). 

\subsection{The structure of the limit cohomology}\label{limitcoh}
By general results, it is known that there exists a mixed Hodge structure $\sH^2_{\lim}$ (defined up to a nilpotent orbit) that represents the limit cohomology of a semistable degeneration of $K3$ surfaces. $\sH^2_{\lim}$ can be defined over $\bZ$ with the same structure as the cohomology of a smooth fiber $\sH^2(X_s,\bZ)$. Thus, $\left(\sH^2_{\lim}\right)_\bZ\cong \Lambda$. For notational simplicity, we use $\left(\sH^2_{\lim}\right)_\bZ$ and $\Lambda$ interchangeably and typically drop the subscript $\bZ$.

\smallskip

The key ingredient associated to a degeneration is the monodromy action, which, in our context, we regard as a unipotent isometry $T\in O(\Lambda)$ acting on $\sH^2_{\lim}$. Then, the logarithm of the monodromy
$$N=\log T=(T-1)-\frac{1}{2}(T-1)^2$$
defines the weight filtration $W_k$ on $(\sH^2_{\lim})_\bQ$. In the case of Type III degenerations of $K3$ surfaces, we have  
 $N^3=0$ and $N^2\neq 0$, and the nontrivial terms of the weight filtration are given by: 
\begin{itemize}
\item[i)] $W_0=\im N^2=\im N\cap \ker N$; 
\item[ii)] $W_2=\ker N^2=\im N+\ker N=W_0^{\perp}$;
\item[iii)] $W_4=(\sH^2_{\lim})_\bQ$
\end{itemize}
In the case of  $K3$ surfaces, $N$ and subspace $W_k$  can be defined over $\bZ$ (cf. \cite[pg. 6]{friedmanscattone}). However, some of the  formulas mentioned above are no longer valid without the use of the saturation. The problem is that $\im N$ and $\im N^2$ are not necessarily primitive sublattices of $\Lambda$. 
In particular, we note that $W_0$ is defined as the saturation of $\im N^2$ (or equivalently the radical of $\ker N$), and thus is a primitive isotropic rank $1$ sublattice of $\Lambda\cong \sH^2_{\lim}$. Up to the action of $O(\Lambda)$, there is only one  choice for $W_0$. It is also easily seen that $\Gr_2^W \sH^2_{\lim}:=W_2/W_0$ is isometric to $M=E_8^{\oplus 2}\oplus U^2$. 

\smallskip

In addition to the weight filtration,  $N$ defines the intermediary subspace $\im N$ and $\ker N$, sitting  between $W_0$ and $W_2$. One has  that  $\im N/W_0$ and $\ker N/W_0$ are mutually orthogonal sublattices in $\Gr_2^W \sH^2_{\lim}$, with $\im N/W_0$ positive definite of rank $1$ and $\ker N/W_0$ indefinite of signature $(1,18)$. Regarded as defined over $\bZ$, $\ker N/W_0$ is automatically primitive in $\Gr_2^W \sH^2_{\lim}$; the same need not be true for $\im N/W_0$. This is essentially due to the fact that the effect of a base change of order $k$ is to replace $T$ by $T^k$ and thus multiply $N$ by $k$ (view $N$ as a matrix).

\smallskip

Friedman--Scattone \cite[Thm. 0.5]{friedmanscattone} have classified the possibilities for $N$ (up to the action of $O(\Lambda)$).  Namely, $N$ is completely specified by two positive integers $k$ and $t$, that can take any values that satisfy the condition $2k^2\mid t$. The integer $k$ is simply the primitivity index of the integral matrix giving the endomorphism $N$ (i.e. $N=k\cdot N'$ with $N'$ integral and primitive).
Geometrically, the index $k$ can be understood as showing that the family $\calX\to \Delta$ was obtained via a base change of order $k$ from an already semistable family (see \cite[Thm. 0.6, (3)]{friedmanscattone} for a precise statement). To define $t$, we let  $\gamma$ be a primitive generator of $W_0\cong \bZ$ and chose $\gamma'\in \Lambda$ such that $\gamma.\gamma'=1$ (N.B. $\gamma'$ is well defined only modulo $W_2$; $\gamma'$ can be regarded as a generator of $(\Gr_4^W\sH^2_{\lim})_{\bZ}\cong \bZ$).  Then, the invariant $t$  of $N$ is defined by 
\begin{equation}
t=\delta.\delta, \textrm{ where } \delta=N \gamma'\in \im N.
\end{equation}
Note that $\delta$ is well defined only modulo $\gamma$, and can be regarded as a generator of $\im N/W_0$. Since $W_0$ is isotropic, $t$ is indeed an invariant of $N$. Geometrically,  $t$  is simply the number of triple points of the central fiber $X_0$ of the degeneration (cf. \cite[Prop. 7.1]{friedmanscattone}).

\smallskip

We close by noting the following explicit formula for $N$:
\begin{equation}
N(x)=(x.\gamma)\delta-(x.\delta)\gamma \textrm{ for all } x\in \Lambda
\end{equation}
(cf. \cite[(1.1)]{friedmanscattone}). In particular, the natural isomorphism of Hodge structures $N^2$ (defined over $\bQ$) is given by 
\begin{eqnarray*}
 N^2: \Gr_4^W\sH^2_{\lim}&\to& \Gr^W_0\sH^2_{\lim}\ ;\\
 \gamma'&\to&-t\cdot \gamma
 \end{eqnarray*}
 fact related to the Picard--Lefschetz transformations (see  \cite[1.8--11]{friedmanscattone}).

\smallskip
 
For further reference, we  summarize the above discussion as follows:

\begin{proposition}[Friedman--Scattone \cite{friedmanscattone}]\label{propfs}
With notations as above, $\im N/W_0$ and $\ker N/W_0$ are mutually orthogonal sublattices in $M=\left(\Gr_2^W\sH^2
_{\lim}\right)_\bZ\cong E_8^{\oplus 2}\oplus U^{\oplus 2}$. The rank $1$ lattice $\im N/W_0$ is spanned by an element $\delta$ of square $t$ and primitivity index $k$ in $M$, where $t$ is the number of triple points of $X_0$. The lattice $\ker N/W_0$ is primitively embedded in $M$ and isometric to $E_8^{\oplus 2}\oplus U \oplus \left(-\frac{t}{k}\right)$. \qed
\end{proposition}

\subsection{Clemens--Schmid exact sequence} In sections \ref{primcoh} and \ref{limitcoh}, we have discussed the cohomology $\sH^2(X_0)$ of the central fiber and respectively the limit cohomology $\sH^2_{\lim}$ of a Type III degeneration of $K3$ surfaces. The  two are closely related by means of the Clemens-Schmid exact sequence (an exact sequence of Hodge structures), which reads:
$$0\to\sH^0_{\lim}\to \sH_4(X_0)\to \sH^2(X_0)\xrightarrow{\mathrm{sp}_2} \sH^2_{\lim}\xrightarrow{N}\sH^2_{\lim}\xrightarrow{\beta}\sH_2(X_0)\to \sH^4(X_0)\to \sH^4_{\lim}\to 0$$ 
In particular, at the level of graded pieces, one obtains the identifications:
$$\overline{\s}^{(0)}_2:\Gr_0^W\sH^2(X_0)\cong\Gr^W_0\sH^2_{\lim}$$
in weight $0$, and respectively
\begin{equation}
\overline{\s}^{(2)}_2:\Gr_2^W\sH^2(X_0)/\Gr_2^W\left(\im\left(\sH_4(X_0)\to \sH^2(X_0)\right)\right)\to\ker N/W_0,
\end{equation}  
in weight $2$. We also recall that in \S\ref{sectcoho} we have identified $\Gr_0^W\sH^2(X_0)=\sH^2(S^2)$ and $\Gr_2^W\sH^2(X_0)=L$.

\smallskip

In general, the above identifications are  valid only over $\bQ$.  In the case of $K3$ surface, as explained in \cite[pg. 24-25]{friedmanscattone}, %(see also \cite[Lemma 3.5]{friedmanannals}), 
there is essentially no change in working over $\bZ$. Specifically, we have the exact sequence of $\bZ$-modules:
\begin{equation}\label{csz}
0\to\bZ\to \sH_4(X_0,\bZ)\to \sH^2(X_0,\bZ)\to \left(\sH^2_{\lim}\right)_\bZ\xrightarrow{N} \left(\sH^2_{\lim}\right)_\bZ 
\end{equation}
(cf. \cite[(4.13)]{friedmanscattone}). Furthermore, $\sH_4(X_0,\bZ)$ is the free module spanned by the fundamental classes $[V_i]$ of the components of $X_0$. The generators $[V_i]$ map to the elements $\xi_i\in L$ via the composition $\sH_4(X_0,\bZ)\to \sH^2(X_0,\bZ)\twoheadrightarrow  \left(\Gr^W_2\sH^2\right)_\bZ=L$.  

\begin{corollary}\label{corfs}
With notations as above, there exist a natural identification between  the lattice $\overline{L}$ associated to $X_0$ and  $\ker N/W_0$. In particular, $\overline{L}$ embeds primitively into $M\cong E_8^{\oplus 2}\oplus U^{\oplus 2}$ as the orthogonal complement of an element $\delta$ of square $t$, equal to the number of triple points of $X_0$. The only extra arithmetic invariant (w.r.t. $\Gamma=O_{-}(M)$) of  the embedding $L\hookrightarrow M$ is the  primitivity index $k$ of $\delta$.
\end{corollary}
\begin{proof}
Since the $0$ weight part of $\sH^2(X_0)$ maps to the $0$ weight part of $\sH^2_{\lim}$ (which is primitive), we obtain map $\phi:L\to \ker N/W_0$ (see the exact sequences (\ref{computeh2}) and (\ref{csz})). Clearly, $\phi$ is surjective and $K\subseteq \ker(\phi)$, with equality when tensoring by $\bQ$. Since $\phi$ respects the intersection products (see \ref{remint}), we get $\ker(\phi)=\Sat(K)$ and the identification follows. The remaining part is Prop. \ref{propfs}.
\end{proof}

\begin{remark}\label{remprimitivity}
By results of Friedman \cite[Thm. 3.4]{friedmanaut}, in the case that $X_0$ is in minus-one-form (see Def. \ref{defminus1}), the primitivity index $k$ has a very simple combinatorial meaning. It measures if the triangulation  $T$ given by dual graph of $X_0$ is obtained as a refinement of another triangulation $T'$ by subdividing each edge $T'$ in $k$ equal intervals and considering the obvious refined triangulation (see  \cite[Def. 3.3]{friedmanaut}).   
\end{remark}

%%%%%%%%%%%%%%%%%%%%%%%%%%%%%%%%%%%%%%%%%%%%%%%%%%%%%%%%%%%%%%
%%% S5
\section{Triangulations of $S^2$ and $K3$ surfaces in minus-one-form}\label{sectminus1}
In the previous sections, we have associated to a Type III $K3$ surface $X_0$  a lattice $\overline{L}$, that reflects the combinatorics of $X_0$. However, the only information contained in $\overline{L}$ is the number of triple points of $X_0$ (see Cor. \ref{corfs}). In this section, we take a closer look at the structure of $\overline{L}$ in the most symmetric case, the case when $X_0$ is in minus-one-form (see Def. \ref{defminus1}). We show that the structure of $\overline{L}$ can be naturally enriched such that it reflects more faithfully the combinatorics of $X_0$. 

\subsection{The canonical Type III $K3$ surface associated to a triangulation}\label{sectcan}
One class of Type III $K3$ surfaces that were  of particular interest in the study of degenerations of $K3$ surfaces are those in minus-one-form. On one hand, the minus-one-forms are general enough: a theorem of Miranda--Morrison \cite{mm} says that any Type III degeneration $K3$ surface can be assumed, after elementary birational modifications, to have the central fiber in minus-one-form.  On the other hand, the minus-one-forms are the simplest from a combinatorial point of view: the only discrete invariant associated to a  surface in minus-one-form is the dual graph. 

\begin{definition}\label{defminus1} Let $X_0$ be a Type III $K3$ surface. We say that $X_0$ is in {\it minus-one-form} if  $D_{ij}^2=-1$ for all double curves, with the exception of the case when $D_{ij}$ is nodal on $\widetilde{V}_i$, when we require $D^2_{ij}=1$. \end{definition}

For minus-one-forms, the additional combinatorial data of $X_0$ given by the self-intersection numbers $D_{ij}^2$ is lost.  To conclude that the combinatorics of $X_0$ is completely determined by the dual graph, we note that the  deformation types of  the components of $X_0$ are determined by the following classification:
 \begin{proposition}\label{structureminus1}
 Let $(\widetilde{V}_i,D_i)$ be a rational surface with an anticanonical cycle of length $d_i$. Assume that $D_{ij}^2=-1$ for all $j$ (or $D_{ij}^2=1$ if $d_i=1$). Then $d_i\le 6$ and $\widetilde{V}_i$ is a del Pezzo surface of degree $d_i$. Furthermore, the orthogonal complement $R_i$ to the sublattice spanned  by the double curves $D_{ij}$ in $\sH^2(\widetilde{V}_i,\bZ)$ is a root lattice of type $E_8$, $E_6$, $D_4$, or $A_2$ for $d_i=1,\dots,4$ respectively (and empty for $d_i=5$ or $6$).
 \end{proposition}
 \begin{proof}
The condition $d_i\le 6$ follows from the Hodge index theorem. It is easy to see that $D_i\in|-K_{\widetilde{V}_i}|$ is big and nef. Thus, $\widetilde{V}_i$ is a (generalized) del Pezzo surface of degree $d_i=D_i^2$. The statement about the orthogonal complement is then standard (see also  Looijenga \cite[Thm. 1.1]{looijengarational} and \cite[\S2]{looijengarational}).
 \end{proof}

In particular, note that the triangulation $T$ associated to a Type III $K3$ surface in minus-one-form is of non-negative combinatorial curvature. Conversely, it is clear that any triangulation $T$ of non-negative curvature can be obtained as the dual graph of   some surface $X_0$ in minus-one-form. This is essentially a $1$-to-$1$ correspondence: the only freedom in choosing $X_0$ to represent a triangulation $T$ is given by the moduli of the components and the choice of the gluing data (i.e. $X_0$ is determined by $T$ up to locally trivial deformations). 

\smallskip

Our first step in enriching the structure of $\overline{L}$ is to note that a surface $X_0$ in minus-one-form is canonically polarized. Specifically, this means  that there exists a polarization $h$ for $X_0$ that is intrinsically determined by the triangulation $T$.  It is the ``anticanonical polarization'' $h$ given by  
\begin{equation}\label{anticanonical}
h:=\sum_i\left(-K_{\widetilde{V_i}}\right)=\sum_i D_i=\sum_{i,j} D_{ij}.
\end{equation}
We note that  
$$h.D_{ij}=D_i.D_{ij}=1=h.D_{ji}.$$ 
Thus, indeed $h\in L$ and we can view the class of $h$ in $\overline{L}$ as a polarization for $X_0$ (see Def. \ref{defpol}). In fact, the condition that the anticanonical class given by (\ref{anticanonical}) is a (combinatorial) Cartier divisor is equivalent to saying that $X_0$ is in minus-one-form. Note also that the degree of $h$ is $3t$, where $t$ is the number of triple points of $X_0$:
$$h^2=\sum_i (D_i)^2=\sum_i d_i=2e =3t.$$

We conclude that to a triangulation $T$ of non-negative combinatorial curvature, we can associate in a natural way a polarized $K3$ surface $(X_0,h)$ in minus-one-form.
\begin{cordef}\label{defk3}
Let $T$ be a triangulation of $S^2$ of non-negative combinatorial curvature. Then, there exists a d-semistable polarized $K3$ surface $(X_0,h)$ in minus-one-form such that:
\begin{itemize}
\item[i)] the dual graph of $X_0$ is isomorphic to $T$;
\item[ii)] $h$ is the class of a semi-ample Cartier divisor on $X_0$;
\item[ii)] the degree of $h$ is $3t$, where $t$ is the number of triangles of $T$; 
\item[iii)]  the degree of each double curve is $1$.
\end{itemize} 
Furthermore, $(X_0,h)$ is uniquely determined by these conditions up to locally trivial\footnote{Strictly speaking, we require local triviality only near the double curves. In other words, we allow rational double points away from the double locus.} deformations that preserve the polarization class $h$. We call $(X_0,h)$ the canonical Type III $K3$ surface associated to $T$. 
\end{cordef}
\begin{proof}
As noted above the condition i) determines the deformation type of $X_0$. Via a Hodge index argument, it is easy to see that the class $h\in \overline{L}$ is uniquely determined by the given numerical conditions. The analytic conditions of d-semistability and Cartier divisor can be always achieved. Finally, the deformation statement is standard (e.g. see \cite{friedmansmooth,friedmanscattone}).
\end{proof}
\begin{remark}
In fact, there is a truly canonical choice for $(X_0,h)$. Namely, the maximal algebraic surface  $X_0^{alg}$ (i.e. $\ker \psi=L$ in the notation of \S\ref{sectcoho}) satisfying the conditions of \ref{defk3} is completely determined by $T$. For the purposes of this paper, one can ignore the analytic structure. Thus, any locally trivial deformation of $X_0^{alg}$ is equally good for us.
\end{remark}

By considering a  polarized Type III $K3$ surfaces $(X_0,h)$, we naturally obtain a tower of lattice embeddings:
\begin{equation}\label{tower}
P\hookrightarrow \overline{L}\hookrightarrow M,
\end{equation}
where $M$ and $\overline{L}$ are as before,  and $P$ is primitive part of $\overline{L}$, i.e. $P:=\langle h\rangle^\perp_{\overline{L}}$. In particular, if $(X_0,h)$ is the canonical surface associated to a triangulation $T$, then the  lattices and the embeddings of (\ref{tower}) are intrinsically determined (up to isometries) by $T$. Thus,  (\ref{tower}) can be regarded as an arithmetic invariant of the triangulation.

\begin{notation} In what follows, $T$ will be a triangulation of non-negative combinatorial curvature, $(X_0,h)$ will be the associated Type III $K3$ surface in minus-one-form, and  $P$, $\overline{L}$ and $M$ will have the same meaning as in this section.  For reader convenience, we recall that $M$ is the fixed lattice $E_8^{\oplus 2}\oplus U^{\oplus 2}$ and $\overline{L}$ is an even lattice of signature $(1,18)$ (see \ref{corfs} for its structure). Also, $\overline{L}$ contains a distinguished element $h$ of norm $3t$ with orthogonal complement $P$. In particular, $P$ is negative definite of rank $18$ and primitively embedded in $\overline{L}$ (and $M$).
\end{notation}

\subsection{The vanishing lattice of a triangulation}\label{sectdegtype} In addition to the number of triangles, the most basic combinatorial invariant associated to a triangulation $T$ of $S^2$ is its {\it degree type}. Namely, if $(d_i)_{i=1,\dots,n}$ denotes the degrees of vertices of $T$, by the Euler formula, we have  
\begin{equation}\label{eulereq}
\sum_i (6-d_i)=12.
\end{equation}
Requiring that the triangulation $T$ is non-negatively curved, i.e. $d_i\le 6$ for all $i$, we obtain $47$ combinatorially distinct solutions to the equation (\ref{eulereq}), which we call the degree type of triangulation. For consistence with Thurston \cite{thurston}, we define:
\begin{definition}
Let $(k_1,\dots, k_l)$ be a partition of $12$ with $1\le k_i\le 5$. We say that a triangulation $T$ of non-negative combinatorial curvature has degree type $(k_1,\dots, k_l)$ if there exists $l$ distinct vertices $v_1,\dots,v_l$ of $T$ with $\deg(v_i)=6-k_i$. If the total number of vertices is $n$, we also say that $T$ has type $(n;k_1,\dots,k_l)$.
\end{definition}

We now note that the degree type of a triangulation $T$ is naturally reflected in the structure of the lattice $P$. Namely, for each component $V_i$ of the associated Type III $K3$ surface $X_0$, we define:
\begin{equation}\label{defri}
R_i:=\left(\Span(\{D_{ij}\}_j)\right)^\perp_{\sH^2(\widetilde{V}_i)}
\end{equation}
and then set
$$R:=\oplus_{i=1}^n R_i$$
(N.B. $R$ is essentially the cohomology supported on the smooth locus of $X_0$). By construction, $R$ is a sublattice of $\langle h\rangle^\perp_L$ with $R\cap \Sat(K)=0$. Thus, $R$ maps isometrically to its image (still denoted by $R$) in  $P\subset \overline{L}=L/\Sat(K)$. On the other hand, by Prop. \ref{structureminus1}, $R_i$  is either empty,  for $d_i=5,6$, or  isometric to a root lattice of type $A_2$, $D_4$, $E_6$, or $E_8$ for $d_i=4,3,2,1$ respectively. We conclude that $R$ is a  root sublattice of $P$; $R$ is  determined by  the degree type of the triangulation $T$ (e.g. $R\cong E_6^{\oplus 3}$ is equivalent to  degree type  $(4,4,4)$ for $T$). In particular, we note that this fact implies that the lattice $P$ contains more combinatorial information about the triangulation $T$ than $\overline{L}$.
 
 \begin{example}\label{exdeg}
 The smallest number of triangles in a triangulation of $S^2$ is $2$. There are two such triangulations, say $T_1$ an $T_2$ of degree type $(5,5,2)$ and $(4,4,4)$ respectively. In both cases the associated lattice $\overline{L}$ is isometric to $E_8^{\oplus 2}\oplus U\oplus (-2)$, but the associated primitive lattices $P_1$ and $P_2$ are not isometric. Specifically, the  root lattices associated to the two triangulations are $R_1\cong E_8^{\oplus 2}\oplus A_2$ and $R_2\cong E_6^{\oplus 3}$ respectively. Since $R_i\subseteq P_i$ and $\rank R_i=\rank P_i=18$, it is easy to see that $P_1=R_1$. Similarly, $P_2$ is an index $3$ overlattice of $R_2$ with the property that the lattice spanned by its roots is  $R_2$.  It follows that $P_1$ is not isometric to $P_2$.
 \end{example}
 
At the same time, the isometry class of $P$  is not enough to determine the triangulation; $P$ reflects only the shape of the triangulation in the sense of Thurston. 
 \begin{example}
For triangulations $T$ of type $(5,5,2)$, the lattice $P$ is determined as in example \ref{exdeg}, i.e. $P\cong E_8^{\oplus 2}\oplus A_2$ and it does not depend on $n$, the number of vertices of the triangulation. The extra data needed to recover the triangulation is the position of $h$ in the lattice $P^\perp_M\cong A_2(1)$. The similar example of degree type $(4,4,4)$ is discussed in detail in Thurston \cite[pg. 518]{thurston}.
\end{example}

\begin{definition}
Let $T$ be a triangulation of non-negative combinatorial curvature of degree type $(d_1,\dots,d_k)$. We call  {\it the vanishing lattice of $T$} the root lattice $R$ obtained as the direct sum of summands of type $A_2$, $D_4$, $E_6$, or $E_8$, one summand for each $d_i=1,\dots,4$ respectively. 
\end{definition}

As noted above $R$ is naturally embedded in $P$ and the name is justified by the following observation.  In \ref{defk3}, we can rigidify the Type III $K3$ surface $X_0$ associated to a triangulation  $T$ by requiring that the components of $X_0$ have no moduli. This can be done,  by imposing $X_0$ to have the maximal configuration of rational double points (see \ref{singminus1}). In this set-up, $R$ is indeed the vanishing cohomology of $X_0$, and should be interpreted as the part of $\overline{L}$ (or $P$) that does not come from $T$.
  
 \begin{lemma}\label{singminus1}
For $d\in\{1,\dots,4\}$, there exists  an anticanonical pair  $(V,D)$ with $D$ having $d$ components, each of self intersection $-1$,  such that $V$ has a unique singular point of type $E_8$, $E_6$, $D_4$ or $A_2$ for $d=1,\dots,4$ respectively. Furthermore, $(V,D)$ is unique up to isomorphism.
\end{lemma}
\begin{proof}
This is a standard application of the theory of del Pezzo surfaces. For example the pair $(V,D)$ with $e=3$ can be constructed by considering in  $\bP^2$ a triangle tangent to a conic. By blowing up twice the tangency points one obtains a degree $3$ del Pezzo with a $D_4$ configuration of nodal curves.
\end{proof}

\subsection{The structure of $\overline{L}$ in the minus-one-case}\label{sectstruclbar}
By definition, the elements of $L$ are combinatorial Cartier divisor. Thus,  for each element $\ell\in L$ there is a well defined degree on the double curves, giving a well defined map:
\begin{eqnarray*}
L&\to&(C^1)^*\cong C_1\\
\ell&\to&\left(d_{ij}\right)_{ij},
\end{eqnarray*}
where $d_{ij}=\deg(\ell_{\mid D_{ij}})$ and $D_{ij}$ is  a double curve.  By Prop. \ref{structureminus1}, the components of $X_0$ are del Pezzo surface of degree at most $6$.  In the case that the degree of the component $V_i$ is at most $5$, the classes of double curves $D_{ij}$ in $\sH^2(\widetilde{V}_i,\bZ)$ are linearly independent and consequently there is no linear relation imposed on the degrees. On the other hand, if the degree is $6$,  it is immediate to see that we get the following linear relations:
\begin{equation}\label{eqhex}
d_{ij_1}-d_{ij_4}=d_{ij_3}-d_{ij_6}=d_{ij_5}-d_{ij_2}
\end{equation}
where $d_{ij_1},\dots,d_{ij_6}$ are the degrees of $\ell$ on the double curves of $V_i$ considered in a consecutive order (w.r.t. to some orientation). We call the equations (\ref{eqhex}) {\it the hexagonal relations} of $X_0$ and conclude that the image of the degree map $L\to C_1$ belongs to the linear subspace: 
$$C:=\ker(C_1\to \bZ^{2h})$$  
where $h$ is the number of hexagonal components (i.e. components of degree $6$) and the map 
$$\psi: C_1\to \bZ^{2h}$$ is given by the equations (\ref{eqhex}).  Clearly, $\psi$ depends only on the combinatorics of the associated triangulation $T$. The fact that $L$ maps to $C=\ker \psi$ should be understood as  a combinatorial analogue to saying that the cone metric on $S^2$ associated  to $T$ (see \S\ref{sectdm}) is flat at the hexagonal points (i.e. the degree $6$ vertices). 

\smallskip

In conclusion, the structure of $L$ is as follows:

\begin{lemma}\label{lemmalc}
With notations as above, we have:
\begin{enumerate}
\item[1)] $\ker(L\to C_1)=R$, where $R$ is the vanishing lattice associated to the triangulation $T$; 
\item[2)] $\im(L\to C_1)=C$.
\end{enumerate}
Thus, $L$ can be viewed as an extension of $C$ by $R$:
\begin{equation}\label{eql0}
0\to R\to L \to C\to 0 
\end{equation}
\end{lemma}
\begin{proof} Let $\ell\in L\subset H=\oplus \sH^2(\widetilde{V}_i,\bZ)$ and write $\ell=\sum\ell_i$ for the component-wise decomposition. By definition, $\ell$ is in the kernel if and only if   $\ell.D_{ij}=\ell.D_{ji}=0$ for all $D_{ij}$. It follows $\ell_i.D_{ij}=0$ for all $j$, and thus $\ell_i\in R_i$ giving that $\ell=\sum \ell_i \in R$. 

Given a degree assignment satisfying (\ref{eqhex}) (i.e. an element of $C$), by using Poincare duality on each component $V_i$ of $X_0$, it is easily seen that there exist elements $\ell_i\in \sH^2(V_i,\bZ)$ satisfying $\deg(\ell_{i\mid D_{ij}})=d_{ij}$ for all $j$ (and $i$). It follows that $\ell=\sum \ell_i$ belongs to $L$ and realizes the given degree assignment. 
\end{proof}

We next note that there exists a canonical section of $L\to C$ defined over $\bQ$. Namely, in the proof of \ref{lemmalc}, we can require $\ell_i\in \sH^2(\widetilde{V}_i,\bQ)$ to be a linear combination of the double curves $D_{ij}$. This requirement defines a well defined section 
$$\sigma:C\to L$$ 
with $\sigma(C)\subseteq R^\perp_L$. Note that the $\ell_i$ are defined over $\bZ$ with the exception of the components of degree $2$, $3$ or $4$, in which case we have to allow division by $2$ or $3$. More precisely, if we consider the elements in $L$ that can be written as a linear combination of double curves and map them to $C$, they span a finite index submodule   $C'\subset C$ such that 
$$C/C'\cong A_R,$$
where $A_R$ is the discriminant group for the root lattice $R$ (it is obtained as a direct sum of summands of type $A_{A_2}\cong \bZ/3$, $A_{D_4}\cong \bZ/2\times \bZ/2$ and $A_{E_6}\cong \bZ/3$). We obtain  that $\sigma_{\mid C'}$ is defined over $\bZ$ and $\sigma(C')$ is contained in the orthogonal complement of $R$ in $L$ (in particular, it has a natural lattice structure). Over $\bQ$, $L$ is a direct sum of $R$ and $\sigma(C')$. Over $\bZ$, $R$ and $\sigma(C')$ are glued (along $A_R$) to give $L$. The point of this discussion is that $L$ is obtained (essentially as a direct sum) from the vanishing lattice $R$ of the triangulation and a lattice $C$ whose generators correspond naturally to the edges of the triangulation $T$ (recall $C\subseteq C_1$, cut out by the hexagonal relations (\ref{eqhex})).

\smallskip

To pass from $L$ to $\overline{L}=L/\Sat(K)$, we note that we have a diagram: 
\begin{equation}\label{eqlbar}
\xymatrix{
K\ar@{^(->}[dr]\ar@{^(->}[d]\\
L\ar@{->}[r]&C,
}
\end{equation}
where the embedding $K\to C$ is induced from the combinatorial map  $C_0\to C_1$  given by sending a vertex $v_i$ of $T$ to 
\begin{equation}\label{eqzeta}
\zeta_i=\sum e_{ij}-\sum e_{jk}
\end{equation}
where $e_{ij}$ are the edges of $T$ incident to the vertex $v_i$, while $e_{jk}$ are the non-incident edges that belong to triangles $t_{ijk}$ containing $v_i$ (N.B. compare the definition of $\zeta_i$ with that of $\xi_i$ from \S\ref{sectcoho}). It is immediate to check that indeed
$$K=\im(C_0\to C_1)\subseteq \ker(C_1\to \bZ^{2h})=C.$$ 
Since $R\cap K=0$ and $K$ is totally isotropic, we can simply replace in the discussion of this section $L$ by $\overline{L}$ and $C$ by $\overline{C}=C/K$. 

\smallskip

For further reference, we also note that 
$$\sum \zeta_i=0$$
and that a choice of orientation of $S^2$ allows us to write (uniquely):
$$\zeta_i=\sum (e_{ij}-e_{jk}),$$
where $e_{jk}$ follows $e_{ij}$ with respect to the orientation of the triangle $t_{ijk}$.

%%%%%%%%%%%%%%%%%%%%%%%%%%%%%%%%%%%%%%%%%%%%%%%%%%%%%%%%%%%%%%
\subsection{Pseudo-geodesics and generators of $\overline{L}$}\label{sectpseudo} 
In order to relate our construction with Thurston \cite{thurston}, we need to find some ``geometric generators'' for $\overline{L}$. Let $\{s_1,\dots,s_k\}$ the singular vertices of the triangulation, i.e. vertices for which the degree is at most $5$ [N.B. $k=n-h$ and $k\le 12$]. The geometric generators for $\overline{L}$ that we consider are the  pseudo-geodesics joining two singular points.  
\begin{definition}
A {\it combinatorial  arc} $\gamma$ for a triangulation $T$ is a collection of oriented edges $(e_1,\dots,e_l)$ such that the end-point $v_j$ of $e_j$ is the starting point of $e_{j+1}$ for $j=1,\dots,l-1$.   An arc is {\it straight} if all the intermediary points $v_j$ (for $j=1,\dots,l-1$) have degree $6$ and the edges $e_j$ and $e_{j+1}$ are opposite in the obvious way (i.e. $e_j^*$ and $e_{j+1}^*$ are opposite edges for the dual hexagon). A  {\it pseudo-geodesic} is a straight arc  joining two singular points. 
\end{definition}

Note that a single edge is a straight arc. Thus, if two singular points are adjacent in a triangulation, then the edge joining them is a pseudo-geodesic. More generally, let $v_0$ be a singular point.  Consider an edge $e_1$ starting from $v_0$. If the end-point $v_1$ of $e_1$ is singular, we stop as we found a pseudo-geodesic. Otherwise, there exist a unique choice for $e_2$ (the opposite of $e_1$ with respect to $v_1$). We continue until we arrive to a singular point $v_l$. The resulting arc is a pseudo-geodesics.  It is clear that no loop is possible, and thus the above procedure terminates by reaching a singular point. Conversely any pseudo-geodesic is obtained in this way. For each singular vertex the number of pseudo-geodesics passing through that vertex is equal to the degree of the vertex.  Thus, we obtain $3k-6$ pseudo-geodesics, where $k$ is the number of singular points. This essentially gives a triangulation of $S^2$ with vertices at the singular points of $T$ and edges being  pseudo-geodesics. This pseudo-triangulation is  a combinatorial analogue of the triangulation by geodesics of Thurston (see \cite[Prop. 3.1]{thurston}). The major difference is that we allow intersections at the hexagonal points (i.e. points of degree $6$, or smooth points for the metric).

\smallskip

By definition a pseudo-geodesic $\gamma_{ij}$ can be viewed as an element of $C_1$ (the orientation does not matter). Since the pseudo-geodesics go straight through a hexagonal point, it is easily seen that the hexagonal relations (\ref{eqhex}) are satisfied. Thus, $\gamma_{ij}\in C$.  Using the canonical section $\sigma:C\to L$, we view $\gamma_{ij}$ as an element of $L$. In particular, we can talk about intersection numbers. In general, $\gamma_{ij}\in L$ is only defined over $\bQ$  (i.e. $\gamma_{ij}\in L\otimes_\bZ \bQ$, but note $6\cdot \gamma_{ij}$ is always defined over $\bZ$). In any case, what we need 
later is the following lemma:

\begin{lemma}\label{lemmageo}
The classes of pseudo-geodesics in $\overline{L}$ span  $\overline{L}/R\otimes_\bZ \bQ$. If $\gamma_{ij}\in L$ is a pseudo-geodesic and  $h\in L$ is the canonical polarization, then their intersection number is computed by  
$$\gamma_{ij}.h=2\cdot \textrm{(length of $\gamma_{ij}$)},$$
where the length of the pseudo-geodesic $\gamma_{ij}$ is equal to the number of edges in the arc defining $\gamma_{ij}$.  \end{lemma}
\begin{proof}
The first claim follows easily by a dimension count. Essentially, the linear space $G$ spanned by pseudo-geodesics is $3k-6$ dimensional, and then $\dim_\bQ K\cap G=k-1$. In general, since $h$ is the anticanonical polarization, the intersection number of $h$ with an element  $\ell\in L$ is twice the sum  of the degrees of $\ell$ on the double curves (e.g. see the proof of \ref{lemmaeven}). The claim is clear.
\end{proof}

\begin{remark}\label{remprimitivity2}
Comparing with Remark \ref{remprimitivity}, it is important to note the effect of refining the triangulation by subdividing each edge in $k$ intervals: the pseudo-geodesics remain the same, but their length is multiplied by $k$. Thus viewing $h$ as an element of $\overline{L}^*$, the subdivision has the effect of multiplying $h$ by $k$.
\end{remark}
The intersection number of two pseudo-geodesics $\gamma_{ij}$ and $\gamma_{kl}$ is easily seen to be:
\begin{equation}\label{eqgeo}
\gamma_{ij}.\gamma_{kl}=n_h+n_1+n_2,
\end{equation}
where $n_h$ is hexagonal contribution and $n_1$ and $n_2$ are the possible contributions from the end-points. The {\it hexagonal contribution} is defined to be the number of times the pseudo-geodesics intersect at  the hexagonal points (the self-intersection is counted twice). The contributions $n_1$ and $n_2$ occur only when the pseudo-geodesics $\gamma_{ij}$ and $\gamma_{kl}$ have one or two of the end vertices in common. Say that the pseudo-geodesics have the vertex $s_i$ (corresponding to the non-hexagonal component $V_i$) in common, then the contribution $n_1$ is computed  by inverting the intersection matrix of the double curves on $V_i$. In other words, we are computing the intersection numbers of the duals to double curves in $\Span(\{D_{ij}\}_j)\subset \sH^2(\widetilde{V}_i,\bQ)$. In the self-intersection case, the  contributions  $n_1$ and $n_2$ from the end-points   are  $-1$, $-\frac{1}{3}$, $0$, $\frac{1}{3}$, $1$ for degrees $5,\dots,1$ respectively. This can be written in a uniform way as follows:
\begin{corollary}\label{selfint}
The self intersection number of a pseudo-geodesic $\gamma_{ij}\in L$  is 
$$\gamma_{ij}^2=-\frac{1}{\sqrt{3}}\left(\cot\frac{k_i\cdot\pi}{6}+\cot\frac{k_j\cdot\pi}{6}\right)+n_h$$
where $k_i=6-\deg(V_i)$ and $V_i$ is the component corresponding to $s_i$ (and similarly for $j$). In particular, if the degrees of the end-points are $5$ and there is no self-intersection at a hexagonal point, then $\gamma_{ij}^2=-2$.  \qed
\end{corollary}

Similarly, if two distinct pseudo-geodesics meet in a singular point $s_i$ of degree $5$, the contribution $n_1$ of $s_i$ to (\ref{eqgeo}) is:
\begin{equation}\label{tempint}
n_1=
\left\{\begin{array}{ll}
0 &\textrm{if the angle between $\gamma_{ij}$ and $\gamma_{il}$ at $s_i$  is }\frac{\pi}{3}\\
1&\textrm{it the angle between $\gamma_{ij}$ and $\gamma_{il}$ at $s_i$ is }\frac{2\pi}{3}
\end{array}
\right. .
\end{equation}
Note that one has to take into account the angle at which the two pseudo-geodesics meet at $s_i$, where the angle is defined as follows. By construction two pseudo-geodesics meeting at the singular point $s_i$ correspond to  two different edges $e_1$ and $e_2$ of the triangulation $T$ starting at $s_i$. By declaring every triangle of $T$ to be equilateral (i.e. every angle to be $\frac{\pi}{3}$), we obtain a well defined minimal angle $\alpha\in\{0,\frac{\pi}{3},\frac{2\pi}{3}\}$ between $e_1$ and $e_2$. We define $\alpha$ to be {\it the intersection angle} of the pseudo-geodesics at $s_i$. Note that $\alpha=0$ forces the two pseudo-geodesics to coincide, and we get the self-intersection case considered above. As before, we obtain:

\begin{corollary}\label{eqint}
Let $\gamma_{ij}$ and $\gamma_{il}$ be two pseudo-geodesics having the end-point $s_i$ in common. The contribution of $s_i$ to  $\gamma_{ij}.\gamma_{il}$  (see (\ref{eqgeo})) is given by:
$$n_1=-\frac{1}{\sqrt{3}}\cdot \frac{\cos \left(\alpha+\frac{k_i\pi}{6}\right)}{\sin \frac{k_i\pi}{6}},$$
where as above $k_i=6-\deg(V_i)$ and $\alpha$ is the angle at $s_i$ between $\gamma_{ij}$ and $\gamma_{il}$. In particular, if $\alpha=0$ we recover the formula of \ref{selfint}.
\end{corollary} 
\begin{proof}
The contribution for degree $5$ is given in (\ref{tempint}) above. For degree $4$ the contribution is  $\frac{1}{3}$ and $\frac{2}{3}$ for $\alpha=\frac{\pi}{3}$ and $\frac{2\pi}{3}$ respectively. For  degrees $3$ and $2$ the angle is automatically $\frac{\pi}{3}$ and we get  $\frac{1}{2}$ and $\frac{2}{3}$ respectively.  By a case by case check, we obtain the given uniform formula for $n_1$. 
\end{proof}

\begin{remark}
From the perspective of the lattice $\overline{L}$, the construction of pseudo-geodesics is natural. Essentially, the pseudo-geodesics are the first non-trivial elements of $L$ that one can construct. The intersection numbers of pseudo-geodesics are easily computed, but the resulting formula of Cor. \ref{eqint} seems rather artificial. It becomes natural only a posteriori when compared with \cite[Prop. 4.6]{loijl}.
\end{remark}
%%%%%%%%%%%%%%%%%%%%%%%%%%%%%%%%%%%%%%%%%%%%%%%%%%%%%%%%%%%%%%
%%% S6
\section{The Eisenstein lattice structure}\label{secteisenstein}
The various lattices that we associated to a triangulation $T$ of $S^2$ via Type III $K3$ surfaces are only $\bZ$-lattices. Here we show, that due to the symmetry of Type III $K3$ surface in minus-one-form, there is an additional natural Eisenstein structure. Locally, if one considers only the double curves incident to a fixed triple point, the Eisenstein structure is quite apparent. The main result of the section (Thm. \ref{mainthm1}) is to show that the obvious local  ``rotation around a triple point'' can be globalized to give an order $3$ fixed-point-free isometry $\rho$ of $M$, thus giving an Eisenstein lattice structure for $M$ (see section \ref{defeisenstein}). Our globalization procedure  roughly corresponds to the construction  of  Thurston \cite{thurston} of a cone metric on $S^2$ by imposing that  every triangle in $T$ is a standard equilateral triangle.  

\begin{theorem}\label{mainthm1}
Let $X_0$ be a Type III $K3$ surface in minus-one-form, with $h$ and $P\hookrightarrow \overline{L}\hookrightarrow M$ as before. Then, there exists an order $3$ fixed-point-free isometry $\rho\in \Gamma$ that leaves $P$ invariant. Additionally, $P^\perp_M$ contains a $\rho$-invariant sublattice $\Delta\cong A_2(-\frac{t}{2})$ with the property that $h=\delta_1+2\delta_2$, for two primitive generators $\delta_1$ and $\delta_2$ of $\Delta$. If we regard $M$ as an Eisenstein lattice $M^\calE$ (see \S\ref{defeisenstein}), the entire arithmetic data is determined by a generator $\delta^\calE\in M^\calE$ of $\Delta^\calE$. Finally, the class of $\delta^\calE\in M^\calE$ modulo $\Gamma^\calE$ depends only on the combinatorics of the dual graph $T$ of $X_0$. 
\end{theorem}

The theorem is proved via several intermediary steps. They are put together in a formal proof in \S\ref{sectconclude}. For reader convenience, we mention that the main construction, the globalization procedure, is in \S\ref{sectmu3} (see esp. Prop. \ref{mu3onc}). 

\subsection{The order  $3$ automorphism defined by a triple point}\label{sectmu3} 
Let $T$ be a triangulation of $S^2$ of non-negative combinatorial curvature. We define a $\bZ$-module $H'=(\bZ^3)^t\cong \bZ^{2e}$ by taking a $\bZ^3$ summand for each triangle of $T$ and generators corresponding to the edges.   We note that a choice of orientation for $S^2$ gives a well defined order $3$ automorphism:
\begin{equation}
\rho:H'\to H'
\end{equation}
induced by the rotation of the edges of each triangle of $T$ according to the orientation. With respect to $\rho$, $H'$ contains an invariant part and an anti-invariant part $A:=\ker(\rho^2+\rho+1)$ of complementary ranks. Clearly, we have $A\cong \bZ^{2t}$ and an induced automorphism:
$$\rho:A\to A \textrm{ with } \rho^3=\mathrm{id} \textrm{ and } \rho(x)\neq x \textrm{ for all } x\in A.$$
We will see that $\rho$ induces an isometry for the primitive lattice $P$. 

\smallskip

To relate $\rho$ to the lattice $\overline{L}$, we first note that there is a natural gluing map $g:H'\to C_1$, which sends $ae'+be''$ to $(a+b) e$, where $a,b\in \bZ$ and $e'$ and $e''$ are the two copies of the edge $e$ in $H'$. By restricting $g$ to the subspace $A$ and composing with the projection map $C_1\to \overline{C}_1$, where $\overline{C}_1=C_1/K$  and $K\hookrightarrow C_1$ is given as in \S\ref{sectstruclbar}, we obtain a map 
\begin{equation}
\phi:A\to \overline{C}_1.
\end{equation}
We then note the following key fact: 

\begin{proposition}\label{mu3onc}
With notations as above, $\ker(\phi)$ is $\rho$-invariant subspace of $A$ of rank $2(n-1)$. In particular, $\rho$ induces a fixed-point-free automorphism of $\im(\phi)\subset \overline{C}_1$, a codimension $1$ (over $\bQ$) subspace of $\overline{C}_1$.
\end{proposition}
\begin{proof} By construction, it is natural to view $A$ as the $A_2^{\oplus t}$ lattice, with an $A_2$ summand for each triangle of $T$ and  the roots corresponding to differences of  edges in the same triangle. Consider next the  $k_i$ triangles that share a vertex $v_i$ of $T$ and the corresponding $A_2$ summands. One notes that there are  well defined roots $\alpha_1,\dots,\alpha_{k_i}$, one for each summand, such that we can write $\alpha_j=e_j-\rho(e_j)$ with both $e_j$ and $\rho(e_j)$ being edges containing the vertex $v_i$. Since the definition of $\rho$ is based on the choice of orientation, it follows that 
$$g(\alpha_1+\dots+\alpha_{k_i})=0.$$
It is easily checked that the linear span 
$$K_1:=\Span(\{\alpha_1+\dots+\alpha_{k_i}\}_i)$$ 
is the kernel of $g_{\mid A}$. Note also that  $K_1$ is the image of a map $\alpha:C_0\to A$, given by sending a vertex $v_i$ to the sum $\alpha_1+\dots+\alpha_{k_i}$.

Let $\widetilde{K}=\ker(\phi)=g_{\mid A}^{-1}(K)$. To show that $\widetilde{K}$ is $\rho$-invariant it suffices to prove that $K_2=\rho(K_1)$ is a section of the induced map $\widetilde{K}\to K$, i.e.  we get the following diagram:  
\begin{equation}\label{diagktilde}
\xymatrix{
\widetilde{K}\ar@{^(->}[dr]&K_2\ar@{_(->}[l]\ar@{=}[r]\ar@{^(->}[d]&K\ar@{^(->}[d]\\
K_1\ar@{^(->}[r]\ar@{^(->}[u]&A\ar@{->}[r]^{g}\ar@{->}[rd]^{\phi}&C_1\ar@{->>}[d]\\
&&\overline{C_1}
}
\end{equation}
with $K_2\hookrightarrow A$ mapping isomorphically onto $K$.  The claim follows by recalling that $K\subset C_1$ is the image of a map $C_0\to C_1$ (see (\ref{eqzeta})), which turns out to coincide with the composition: 
$$C_0\xrightarrow{\alpha}A\xrightarrow{\rho}A\xrightarrow{g} C_1.$$
Explicitly, using the earlier notations, if we let $\beta_j=\rho(\alpha_j)$, we have
$$g(\beta_1+\dots+\beta_{k_i})=\zeta_i,$$
where $\zeta_i$ are the generators of $K$ given by (\ref{eqzeta}).
\end{proof}

It is clear that $\widetilde{Q}:=\im(A\to C_1)$ coincides with the kernel of the map $C_1\xrightarrow{(1,\dots,1)} \bZ$. We then let
\begin{equation}
\widehat{Q}:=\widetilde{Q}/K=\im(\phi).
\end{equation}
The previous proposition says  that there is an induced order $3$ automorphism: 
\begin{equation}
\rho:\widehat{Q}\to\widehat{Q} \textrm{ with } \rho^3=\mathrm{id} \textrm{ and } \rho(x)\neq x \textrm{ for all } x\in \widetilde{Q}.
\end{equation}
(by abuse of notation, we use $\rho$ to denote all induced automorphisms from $\rho:H'\to H'$). We are interested in an order $3$ automorphism on the primitive part of the cohomology $P$. Thus we are interested not in $\widehat{Q}$, but rather in 
$$Q:=(C\cap \widetilde{Q})/K,$$ 
or equivalently the kernel of the map $\widehat{Q}\to Z^{2e}$ induced from $C_1\to \bZ^{2e}$ (see \S\ref{sectstruclbar}, esp. (\ref{eqhex})). The following lemma shows that $\rho$ descends to an automorphism of $Q$. In other words, the order $3$ automorphism $\rho$ introduced here is compatible with the flattening of the hexagonal points.

\begin{lemma}\label{mu3onq}
With notations as above,  $Q$ is a $\rho$-invariant subspace of $\widehat{Q}$. In particular, the restriction $\rho_{\mid Q}$ induces a fixed-point-free order $3$ automorphism of $Q$.
\end{lemma}
\begin{proof}
The statement of the lemma is saying that for $x\in Q$, we have $\rho(x)\in Q$. Since $\rho$ is defined by means of $A$, this is equivalent to check that if $\hat{x}$ is any lift of $x$ to $A$, then $\rho(\hat{x})$ is  in the kernel of the composite map $A\xrightarrow{\phi} \widehat{Q}\to \bZ^{2h}$. Concretely, this  means that given a degree assignment  $x=(d_{ij})_{ij}$ that satisfies (\ref{eqhex}) and a lift $\hat{x}=(d_{ij}',d_{ij}'')_{ij}$ to $A$ (thus $d_{ij}=d_{ij}'+d_{ij}''$), the rotated degree assignment $\rho(x)= (\rho(d_{ij}')+\rho(d_{ij}''))_{ij}$ still satisfies (\ref{eqhex}). Clearly, it suffices to check this one hexagonal component $V_i$  at a time. Considering the $6$ triangles corresponding to the $6$ triple points of $V_i$, we get $18$ edges (i.e. generators of $H'$) involved in the computation. Since we are restricted to $A$, there are also $6$ linear relations. It follows that the rotated degree assignment for $V_i$: 
\begin{equation}\label{newdeg}
\left(\rho(d_{ij_1}')+\rho(d_{ij_2}''), \dots, \rho(d_{ij_6}')+\rho(d_{ij_6}'')\right)
\end{equation}
depends only on  $(d_{ij_1}',d_{ij_2}'', \dots, d_{ij_6}'')$, with $d_{ij_k}=d_{ij_k}'+d_{ij_k}''$ for $k=1,\dots,6$. For example, for appropriate numbering, we have:
$$\rho(d_{ij_1}')=-d_{ij_1}'-d_{ij_6}'' \textrm{ and } \rho(d_{ij_1}'')=d_{ij_2}'$$ 
It is then straight-forward to check that $(d_{ij_1}, \dots, d_{ij_6})$ satisfy (\ref{eqhex}) if and only if the  assignment resulting by applying $\rho$ (see (\ref{newdeg})) satisfies  (\ref{eqhex}).
\end{proof}

To conclude, starting from the obvious local rotation of the edges around a triple point, we obtained a fixed point free order $3$ automorphism on the $\bZ$-module $Q$ (where $Q$ coincides with $P$ modulo the vanishing lattice $R$, see \S\ref{sectstruclbar}). This is done by first globalizing (using an orientation of $S^2$) to an automorphism of $H'$ and then descending (with intermediary steps $A$ and $\widehat{Q}$) to $Q$ (cf. \ref{mu3onc} and \ref{mu3onq}).

\subsection{The Eisenstein lattice on the primitive cohomology} In the previous section we have introduced a fixed-point-free order $3$ automorphism $\rho$ on the $\bZ$-module $Q$. By declaring $\rho$ to be the multiplication by $\omega$, $Q$ becomes a module over the Eisenstein integers $\calE=\bZ[\omega]$. We now show that the structure of Eisenstein module can be extended in a compatible way to the primitive lattice $P$. 

\smallskip

We recall that the relation between $P$ and $Q$ is given by the exact sequence 
\begin{equation}\label{extension1}
0\to R\to P\to Q\to 0
\end{equation}
which is simply (\ref{eqlbar}) restricted to the primitive part (see \S\ref{sectstruclbar} and \S\ref{sectdegtype}). Since $P$ is built out of $R$ and $P$ and both have natural Eisenstein structures, we obtain:

\begin{proposition}\label{extendp}
 There exists a fixed-point-free order $3$ isometry $\widetilde{\rho}$ of $P$ (unique up to isometries of $P$) extending the automorphism $\rho$ of $Q$. In particular, $P$ has the structure of an Eisenstein lattice, that we denote by $P^\calE$.
 \end{proposition}
 \begin{proof} As noted in \S\ref{sectmu3}, over $\bQ$, we can naturally identify $Q$ with the orthogonal complement of $R$ in $P$. Thus, $Q$ has an induced negative-definite intersection pairing such that $P$ becomes an orthogonal direct sum of $R$ and $Q$.  Since $R$ and $Q$ can be viewed as Eisenstein modules (cf. Ex. \ref{deflat} and \S\ref{sectmu3} respectively), we obtain that  $P_\bQ:=P\otimes_\bZ\bQ$ has a natural structure of $\bQ[\om]$ vector space. The remaining statement (over $\bQ$) is that this structure is compatible with the intersection pairing on $P$. Clearly, this is equivalent to saying that $\rho$ is an isometry on $Q$ with respect to the intersection pairing $\langle\cdot,\cdot\rangle$ induced from $P$. On the other hand, by the construction of \S\ref{sectmu3}, $Q$ is obtained from $A$ (essentially as a quotient). Considering $A$ with the obvious lattice structure (isometric to $A_2^{\oplus t}$, see  the proof of \ref{mu3onc}), we get  that $Q$ can be endowed with a second negative-definite intersection pairing $\langle\cdot,\cdot\rangle'$, this time inherited from $A$. By construction, $\rho$ is an isometry of $A$.  It follows that $\rho$ is an isometry of $Q$ with respect to $\langle\cdot,\cdot\rangle'$. We conclude that the statement about the lifting of $\rho$ to an isometry $\widetilde{\rho}$ of $P$ would follow provided that we show that $\langle\cdot,\cdot\rangle$ and $\langle\cdot,\cdot\rangle'$ coincide up to scaling. For clarity, we prove this only in the case when there are no vertices of degree $6$ (the general case is handled along the lines of \S\ref{sectpseudo}). 

Assuming that there are no hexagonal components for $X_0$, the pseudo-geodesiscs of \S\ref{sectpseudo} coincide with the edges $e_{ij}$ of the triangulation $T$. Clearly, $Q\subseteq P$ is generated by elements of type $x_{ijk}:=e_{ij}-e_{ik}$, with $e_{ij}$ and $e_{ik}$ being $2$ incident edges of a triangle $t_{ijk}$. The various intersection numbers $\langle x_{ijk},x_{i'j'k'}\rangle$ are then computed by the formulas of \S\ref{sectpseudo} (esp. Cor. \ref{eqint}). On the other hand, we recall that $A\twoheadrightarrow Q$ with kernel $\widetilde{K}$ (in the case of no hexagonal component). Then, $\langle x_{ijk},x_{i'j'k'}\rangle'$.  Thus, $\langle x_{ijk},x_{i'j'k'}\rangle'$ is computed by picking-up lifts $y_{ijk}\in \widetilde{K}^\perp_A$ of  $x_{ijk}$. The computation can be done one triangle at a time. A case by case check gives the needed statement.
 
 By the above discussion, we can construct over $\bQ$ an isometry $\widetilde{\rho}$ extending $\rho$. The remaining step is to argue that $\rho$ can be defined over $\bZ$. For this, recall that the section $\sigma:C\to L$ of \S\ref{sectstruclbar} can be defined over $\bZ$ when we restrict to the submodule $C'\subseteq C$ with $C/C'\cong A_R$. It is not hard to see that $\rho$ induces a $\mu_3$-action on $C/C'$. The statement about the $\widetilde{\rho}$ being defined over $\bZ$ is equivalent to saying that this $\mu_3$-action on $C/C'\cong A_R$ can be identified with the $\mu_3$-action on $A_R$ induced from Eisenstein structure on $R$. This is checked one special component (i.e. of degree $2$, $3$, or $4$) at a time.
 \end{proof}

 \subsection{Taking into account the polarization} The orthogonal complement of $P$ in $M$ is positive definite of rank $2$. We note here that it can be endowed with an Eisenstein structure.
\begin{lemma}\label{lemmaa2}
Let $(X_0,h)$ and $\overline{L}\hookrightarrow M$ be the objects associated to a triangulation as before. Then, there exists an element $\delta\in \overline{L}^\perp_M$ with $\delta^2=t$ and such that $h-\delta$ is $2$-divisible. Thus, if  $\delta'=\frac{h-\delta}{2}$, then $\Delta:=\Span(\delta,\delta')$ is isometric to $A_2(-\frac{t}{2})$.
\end{lemma}
\begin{proof}
By Cor. \ref{corfs}, we know that $\overline{L}^\perp_M$ contains an element $\delta$ (unique up to $\pm$) with  $\delta^2=t$. Since $M$ is unimodular, we have 
\begin{equation}\label{isomdisc}
A_{\overline{L}}\cong A_{\Sat(\langle \delta\rangle)}\cong \bZ/t',
\end{equation}
where $t'=t/k$, with $k$ being the primitivity index of $\delta$. Clearly, the class of $\frac{\delta}{2}$ is an element of order $2$ in $A_{\Sat(\langle \delta\rangle)}$. On the other hand, we noted (see \ref{lemmaeven}) that $h.\overline{L}\equiv 0 \mod 2$. It follows that $\frac{h}{2}\in \overline{L}^*$, giving an element of order $2$ of $A_{\overline{L}}\cong \overline{L}^*/\overline{L}$.  The element $h-\delta$ is divisible by $2$ in $M$ iff the two order $2$ elements in $A_{\overline{L}}$
are identified via the isomorphism (\ref{isomdisc}). Since the  discriminant group is cyclic, this follows provided that they are either simultaneously trivial or non-trivial. In general, we get non-triviality, i.e. $\delta$ and $h$ are not divisible by $2$ (e.g. it is always true for trivial reasons if $t\not\equiv 0 \mod 8$). The non-primitive case needs special handling, but this can be done by following the arguments of Friedman \cite{friedmanaut} on the existence of special bands of hexagons (see Remarks \ref{remprimitivity} and \ref{remprimitivity2}).
\end{proof}

On $\Delta\cong A_2(-\frac{t}{2})$  there is an obvious Eisenstein lattice structure (unique up to isometries).  We then have
\begin{equation}
P\oplus \Delta\subseteq M
\end{equation}
is a finite index sublattice. By arguments similar to those of Prop. \ref{extendp}, we get that  $M$ can be given the structure of an Eisenstein lattice such that $P^\calE$ and $\Delta^\calE$ are mutually orthogonal sublattices.

\subsection{Proof of Thm. \ref{mainthm1}}\label{sectconclude}
As described above, the lattice $M$ is obtained by gluing the geometrically meaningful pieces $R$, $Q$ and $\Delta$ [N.B. they are mutually orthogonal sublattices in $M$ with $R\oplus Q\oplus \Delta\subset M$ of finite index]. The Eisenstein structure for $R$ and $\Delta$ is clear due to their root lattice type structure (see \S\ref{sectdegtype}, Lemma \ref{lemmaa2}, and Ex. \ref{deflat}). The only interesting aspect is the Eisenstein structure on $Q$. This is introduced in \S\ref{sectmu3} (esp. Lemma \ref{mu3onc}). Prop. \ref{extendp} then shows that the Eisenstein structure on $Q$ is compatible with the intersection form. By the same proposition, the Eisenstein lattices $Q$ and $R$ glue to give an Eisenstein structure on $P$. A similar argument gives that $P$ and $\Delta$ glue to give an Eisenstein structure for $M$. By the discussion of section \ref{defeisenstein}, we obtain the first part of the theorem.  By Cor. \ref{corunique}, up to isometries, there can be only one Eisenstein lattice on $M$, denoted $M^\calE$ in this paper. Since $P$ is the orthogonal complement of $\Delta$ in $M$ and $\Delta$ is of rank $1$ when regarded as an Eisenstein lattice, the entire arithmetic data is determined by a point:
$$\delta^\calE\in M^\calE,$$ 
where $\delta^\calE$ is a generator of $\Delta^\calE$. By construction,  the norm of $\delta^\calE$ is $t$, the number of triple points of $X_0$.

\smallskip

Clearly, the entire construction depends only on the isotopy class of the triangulation $T$, the dual graph of $X_0$. Namely, 
 in section \ref{sectminus1}, we have noted that the isotopy classes of the triangulations of non-negative combinatorial curvature of $S^2$ correspond in one-to-one way with the (locally trivial) deformation types of polarized  Type III $K3$ surfaces  $(X_0,h)$ in minus-one-form. By considering the cohomology, the pair $(X_0,h)$ determines the tower of lattice embeddings:
\begin{equation}\label{embpm}
P\hookrightarrow \overline{L}\hookrightarrow M.
\end{equation}
Since (\ref{embpm}) is invariant under locally trivial deformations, the embeddings of (\ref{embpm}) are determined by $T$.
Finally, the only essential choice in the construction of the Eisenstein structure is the order $3$ automorphism $\rho$ on $Q$, but this is constructed in terms of the triangulation $T$. \qed

%%%%%%%%%%%%%%%%%%%%%%%%%%%%%%%%%%%%%%%%%%%%%%%%%%%%%%%%%%%%%
 %%% S7
 \section{Equivalence of the two invariants associated to a triangulation}\label{sectequivalence}
 W. Thurston \cite{thurston} associated to a triangulation $T$ of $S^2$ a complete invariant $\delta_1^\calE$ of arithmetic nature (Thm. \ref{thmthurston}). As explained in \S\ref{sectdm}, this is done by using the moduli of $12$ points in $\bP^1$. In this paper, using a different geometric construction, involving degenerate $K3$ surfaces,  we have obtained another arithmetic invariant $\delta_2^\calE$ associated to $T$ (Thm. \ref{mainthm1}). The two invariants $\delta_1^\calE$ and $\delta_2^\calE$ have the same nature; we now show that they actually coincide. 
 
 \begin{theorem}\label{mainthm2}
Let $T$ be a triangulation of $S^2$ of non-negative combinatorial curvature. Consider the elements $\delta_1^\calE$ and $\delta_2^\calE$ in $M^\calE$ associated to $T$ by  the construction via Type III $K3$ surfaces (Thm. \ref{mainthm1}) and respectively Thurston (Thm. \ref{thmthurston}). Then, there exists an isometry  $\phi:M^\calE\to M^\calE$ such that $\phi(\delta^\calE_1)=\delta^\calE_2$.
 \end{theorem}
 
 The lattices occurring in the two constructions are both abstractly isomorphic to $M^\calE$. The idea of the proof is to identify a system of generators for $M^\calE$, thought of as cycles, that can be interpreted geometrically in both constructions. The theorem then follows by showing that the two invariants evaluate to the same thing in the two cases. The generators that we consider are the pseudo-geodesics introduced in \S\ref{sectpseudo}. They can be interpreted as cycles in Thurston construction by using a slight modification of the arguments of Looijenga \cite{loijl} on  Lauricella arc systems. This is done in \S\ref{sectcycles} below. We also check the compatibility of the lattice structures. The main ingredient that makes the comparison possible is the combinatorial nature of the pseudo-geodesics and the fact that the computations on both sides are reduced to basically counting. The proof of the theorem is completed in \S\ref{proof2}.
 
 \begin{remark}
 It is likely that a more geometric proof of Thm. \ref{mainthm2} exists. As we noted in remark \ref{remk3}, a construction equivalent to  that of Deligne--Mostow can be done with elliptic $K3$ surfaces with a $\mu_3$-automorphism. Thus, both sides occurring in Thm. \ref{mainthm2} can be represented by $K3$ surfaces or degenerations of those. Thus, it is likely that the two can be directly related by degenerations arguments. 
 \end{remark}
 
%%%%%%%%%%%%%%%%%%%%%%%%%%%%%%%%%%%%%%%%%%%%%%%%%%%%%%%%%%%%%%
\subsection{Eigencycles and Triangulations}\label{sectcycles}
As mentioned in \S\ref{identn}, the space of cycles occurring in Thurston is  $\sH_1(C,\calE)_{\overline{\zeta}}$, where $C\to \bP^1$ is the $\mu_6$-cover of $\bP^1$ branched at the singular points and  $\zeta=-\overline{\omega}$ is a primitive root of unity of order $6$. Given a triangulation $T$ of $S^2\cong \bP^1$, and the cyclic  $\mu_6$-cover $C\to \bP^1$ (see \S\ref{identn}), we obtain a triangulation $T'$ of $C$. We denote by $C_j'$ the corresponding chains and then 
$$C_j'(\calE):=C_j'\otimes_\bZ(\calE).$$ 
The covering transformation $\psi:C\to C$, given by $\psi(z,t)=(\zeta\cdot z,t)$, induces an action of $\mu_6$ on $C_j'$, that we denote by $r$. Assume first that there is no hexagonal point for the triangulation $T$, it is immediate to see that we get an exact sequence:
 \begin{equation}\label{eigen1}
0\to C_2'(\calE)_{\overline{\zeta}}\to C_1'(\calE)_{\overline{\zeta}}\to \sH^1(C,\calE)_{\overline{\zeta}}\to 0,
\end{equation}
where the subscript $\overline{\zeta}$ denotes the $\overline{\zeta}$-eigenspace [N.B. the orientation class is $\mu_6$-invariant]. In general, one obtains a complex that computes $\sH^1(C,\calE)_{\overline{\zeta}}$. We  get the following diagram:
\begin{equation}\label{eigen2}
\xymatrix{
&0\ar@{->}[d]&0\ar@{->}[d]\\
&C_2'(\calE)_{\overline{\zeta}}\ar@{->}[d]\ar@{=}[r]&C_2'(\calE)_{\overline{\zeta}}\ar@{->}[d]\\
0\ar@{->}[r]&Z_1(\calE)_{\overline{\zeta}}\ar@{->}[r]\ar@{->}[d]&C_1'(\calE)_{\overline{\zeta}}\ar@{->}[r]\ar@{->}[d]&C_0'(\calE)_{\overline{\zeta}}\ar@{->}[r]\ar@{=}[d]&0\\
0\ar@{->}[r]&\sH_1(C,\calE)_{\overline{\zeta}}\ar@{->}[d]\ar@{->}[r]&\widetilde{Q}(\calE)_{\overline{\zeta}}\ar@{->}[d]\ar@{->}[r]&C_0'(\calE)_{\overline{\zeta}}\ar@{->}[r]&0\\
&0&0
}
\end{equation}
where $\widetilde{Q}(\calE)_{\overline{\zeta}}=C_1'(\calE)_{\overline{\zeta}}/C_2'(\calE)_{\overline{\zeta}}$. It is not hard to see (and essentially discussed in \cite[Sect. 2]{delignemostow}) that the dimensions of the spaces involved are:
\begin{eqnarray*}
\dim_\bC C_2'(\bC)_{\overline{\zeta}}&=&t\\
\dim_\bC C_1'(\bC)_{\overline{\zeta}}&=&e\\
\dim_\bC C_0'(\bC)_{\overline{\zeta}}&=&h
\end{eqnarray*}
giving that the dimension of $\sH_1(C,\bC)_{\overline{\zeta}}=\sH_1(C,\calE)_{\overline{\zeta}}\otimes_\calE \bC$ is $(n-h)-2$.

\smallskip

We are interested in relating to our construction, which involves real lattices. To do this, we note that the eigenchains are obtained   via the natural maps:
\begin{equation}
\xymatrix{
C_2'\ar@{->}[r]\ar@{^{(}->}[d]&C_1'\ar@{->}[r]\ar@{^{(}->}[d]&C_0'\ar@{^{(}->}[d]\\
C_2'(\calE)\ar@{->}[r]\ar@{->>}[d]&C_1'(\calE)\ar@{->}[r]\ar@{->>}[d]&C_0'(\calE)\ar@{->>}[d]\\
C_2'(\calE)_{\overline{\zeta}}\ar@{->}[r]&C_1'(\calE)_{\overline{\zeta}}\ar@{->}[r]&C_0'(\calE)_{\overline{\zeta}}
}
\end{equation}
explicitly given by 
\begin{equation}\label{take1}
x\in C_j' \to \widetilde{x}=\pi_{\overline{\zeta}}(x):=\frac{1}{6}\sum_{i=0}^5 \zeta^i \cdot r^i (x)  \in C_j'(\calE)_{\overline{\zeta}},
\end{equation}
where as mentioned above $r=\psi_*:C_j'\to C_j'$ (for $j=0,1,2$). As essentially discussed in \S\ref{defeisenstein}, this gives an identification between 
$$\pi_{\overline{\zeta}}:(C_j')_{x^2-x+1}\cong C_j'(\calE)_{\overline{\zeta}},$$
where
$$(C_j')_{x^2-x+1}=\ker(r^2-r+1:C_j'\to C_j')$$
is defined over $\bZ$. Since $r^2-r+1=0$ on $(C_j')_{x^2-x+1}$, (\ref{take1}) simplifies to 
\begin{equation}\label{take2}
\widetilde{x}=\frac{1}{2}\cdot x +\frac{i}{2\sqrt{3}} \cdot (x+2\cdot r^2(x))
\end{equation}
with inverse $2\re:C_j'(\calE)_{\overline{\zeta}}\to (C_j')_{x^2-x+1}$.

\subsection{Pseudo-geodesics as eigencycles}\label{sectlift} 
We now lift our pseudo-geodesics of \S\ref{sectpseudo} to eigencycles and show that the intersection formulas of \S\ref{sectpseudo} (esp. Cor. \ref{eqint}) are the real part of the intersection formulas obtained by using the Hermitian form on $\sH_1(C,\calE)_{\overline{\zeta}}$. For clarity, we restrict first to the case when there are no hexagonal points. Thus, we assume that the pseudo-geodesics are simply edges of the triangulation (i.e. generators of $C_1$) and they coincide with the geodesics of Thurston \cite[Prop. 3.1]{thurston}.

\smallskip

By the discussion of \S\ref{sectstruclbar}, we know that (over $\bQ$) $C_1=R^\perp_L$ and thus it has a lattice structure. The radical of $C_1$ is $K$, which we view as the image of a map $C_0\to C_1$ defined in \S\ref{sectstruclbar}. Similarly, by extending by $0$ the hermitian form on $\sH_1(C,\calE)_{\overline{\zeta}}$, we obtain a hermitian form $h$ on $C_1'(\calE)_{\overline{\zeta}}$ with radical $C_2'(\calE)_{\overline{\zeta}}$ (see (\ref{eigen1})). We assume that the hermitian form $h$ is scaled as in Cor. \ref{corallcock} (so that it is the Hermitian form used by Thurston \cite{thurston}). The idea of our construction is to lift the edges $e_{ij}$ of $T$ to eigenchains $\widetilde{e_{ij}}$ and to compare the intersections in $C_1$ and $C_1'(\calE)_{\overline{\zeta}}$.
 
 \smallskip
 
 Generally speaking, given an edge $e$ of $T$ we can consider lift $e'$ to $T'$. By the discussion of the previous section, we then obtain an  eigenchain associated to $e$ by:
 \begin{equation}\label{eigenchain}
 \widetilde{e}=\pi_{\overline{\zeta}} (e')=\frac{1}{6}\sum_{i=0}^5 \zeta^i \cdot r^i (x).
 \end{equation}
The problem with this construction is that given $e$, there are $6$ different lifts $e'$. However, they differ  by a power of $r$ (i.e. a second lift $e''$ satisfies $e''=r^i(e')$ for some i). Thus, $\widetilde{e}$ is  well-defined by $e$ up to multiplication by $\zeta^i$ for some $i$. Given two edges $e_1$ and $e_2$ our goal is to find lifts $\widetilde{e_1}$ and $\widetilde{e_2}$ such that 
\begin{equation}\label{req1}
\langle e_1,e_2\rangle=\re\ h(\widetilde{e_1},\widetilde{e_2})
\end{equation}
where on the left side we consider the intersection form on $C_1$ and on the right side the Hermitian form on $C_1'(\calE)_{\overline{\zeta}}$. More precisely, we are looking for a section 
\begin{equation}\label{section}
s:C_1\to C_1'(\calE)_{\overline{\zeta}}
\end{equation}
such that the requirement (\ref{req1}) is satisfied. Let us note first that this makes sense. The intersection number $\langle e_1,e_2\rangle$ was computed in \S\ref{sectpseudo}. On the other hand, we have the following:
\begin{lemma}\label{lemmal}
With assumptions as above, let $e_{ij}$ and $e_{jk}$ be two incident edges making an angle of $\alpha\in \{0,\frac{\pi}{3},\dots,\frac{5\pi}{3}\}$ measured in the direction of the orientation of $T$ (i.e. the arc $(e_{ij},e_{jk})$ used to measure $\alpha$ is correctly oriented, see also \S\ref{sectpseudo}). Then there exists lifts such that the associated eigenchains satisfy
\begin{equation}\label{req2}
h(\widetilde{e_{ij}},\widetilde{e_{jk}})=\frac{1}{\sqrt{3}}\cdot \frac{1}{\sin \frac{k_j \pi}{6}}\cdot e^{(\pi-\frac{k_j \pi}{6} -\alpha) i}.
\end{equation}
In particular, (\ref{req1}) is satisfied for these lifts.
\end{lemma}  
\begin{proof}
With our assumptions we are in the situation discussed by Looijenga \cite[\S4.3]{loijl}. Specifically,  $e_{ij}$ and $e_{jk}$ can be viewed as two consecutive arcs in a Lauricella arc system. \cite[Prop. 5.4]{loijl} computes  the value of the hermitian form for certain lifts of these arcs. Specifically, scaling the Hermitian form $H$ used in \cite{loijl} by $-\frac{4}{\sqrt{3}}$ to obtain  $h$  and using the consecutive arcs $\epsilon'_j$ and $\epsilon'_{j+1}$ passing through the vertex $v_j$ (see the proof of Cor. \ref{corallcock} for further details), we have 
\begin{equation}\label{tempo}
h(\epsilon'_{j},\epsilon'_{j+1})=\frac{1}{\sqrt{3}}\cdot \frac{1}{\sin\frac{k_j\pi}{6}}\cdot e^{-\frac{k_j\pi}{6}i}.
\end{equation}
[N.B. the lifts $\epsilon'_{j}$ and $\epsilon'_{j+1}$ depend on choices. What we need is that there exist some lifts for which (\ref{tempo}) holds and that, due to the continuity property outside the Lauricella arc system, a triangle lifts to a triangle.]

Let $\widehat{e_{ij}}=\epsilon'_{j}$ and $\widehat{e_{jk}}=\epsilon'_{j+1}$ be the two lifts that satisfy (\ref{tempo}). We define the lifts $\widetilde{e_{ij}}$ and $\widetilde{e_{jk}}$ by rotating these natural lifts: 
\begin{eqnarray}
\notag \widetilde{e_{ij}}&=& \widehat{e_{ij}}\\
\widetilde{e_{jk}}&=& e^{-(\pi-\alpha)i}\cdot  \widehat{e_{jk}}\label{scale}
\end{eqnarray}
Note that $e^{-(\pi-\alpha)i}$ is always a unit in $\calE$ (a power of $\zeta$). Roughly speaking, $\widehat{e_{ij}}=\epsilon'_{j}$ and $\widehat{e_{jk}}=\epsilon'_{j+1}$  should be thought as coming from arcs in $\bP^1$ making an angle of $\pi$; the scaling (\ref{scale}) makes the angle between arcs $\alpha$. 

With the definition of the lifts as in (\ref{scale}), we obtain:
\begin{equation*}
h(\widetilde{e_{ij}},\widetilde{e_{jk}})=h(\epsilon_{j},e^{-(\pi-\alpha)i}\epsilon_{j+1})=\frac{1}{\sqrt{3}}\cdot \frac{1}{\sin\frac{k_j\pi}{6}}\cdot e^{(\pi-\frac{k_j\pi}{6}-\alpha)i}.
\end{equation*}
In particular, for the real part we have
$$\langle e_{ij},e_{jk}\rangle=\re\ h(\widetilde{e_{ij}},\widetilde{e_{jk}})=-\frac{1}{\sqrt{3}}\frac{\cos \left(\alpha+\frac{k_j\pi}{6}\right)}{\sin \frac{k_j\pi}{6}}.$$
Note that the computation is valid also in the case that the two edges coincide, i.e. $\alpha=0$, provided that we account for both end-points:
\begin{eqnarray}\label{tempo2}
h(\widetilde{e_{ij}},\widetilde{e_{ij}})&=&\frac{1}{\sqrt{3}}\cdot\left( \frac{1}{\sin \frac{k_i \pi}{6}}\cdot e^{(\pi-\frac{k_i \pi}{6} -\alpha) i}+
 \frac{1}{\sin \frac{k_j \pi}{6}}\cdot e^{(\pi-\frac{k_j \pi}{6} -\alpha) i}\right)\\
&=&-\frac{1}{\sqrt{3}}\cdot\left( \cot  \frac{k_i \pi}{6}+\cot \frac{k_j \pi}{6}\right)\notag
\end{eqnarray}
 (see also Cor. \ref{selfint} and \cite[Prop. 5.4]{loijl}).
\end{proof}

\begin{remark}
We emphasize that to assure consistence, an orientation of $S^2$ is fixed throughout. The angle $\alpha$ is measured with respect to this orientation. Reversing the orientation has the effect of taking the conjugate in (\ref{req2}):
$$\overline{ e^{(\pi-\frac{k_j \pi}{6} -\beta) i}}= e^{(\pi-\frac{k_j \pi}{6} -\alpha) i}\ \textrm{ for }\ \alpha+\beta=2\pi-\frac{k_j \pi}{3}\textrm{ (the cone angle at $v_j$).}$$
In other words, as expected:
\begin{equation*}
h(\widetilde{e_{ij}},\widetilde{e_{jk}})=\overline{h(\widetilde{e_{jk}},\widetilde{e_{ij}})}.
\end{equation*} 
\end{remark}

Now the construction procedure of the sections $s$ of (\ref{section}) satisfying (\ref{req1}) is clear. We start by choosing an arbitrary lift $\widetilde{e_{ij}}$ for a first edge $e_{ij}$. Then for the edges $e_{ik}$ incident to $e_{ij}$ we chose the unique lift $\widetilde{e_{ik}}$ that satisfies Lemma \ref{lemmal}. We continue inductively (in an arbitrary way) until we have chosen lifts for all $e_{ij}$. By extending by linearity, we get a section $s$ as needed. We conclude:

\begin{proposition}
Assuming that there is no point of degree $6$ for the triangulation $T$, there exists a $\bZ$-linear map 
$$s:C_1\to C_1'(\calE)_{\overline{\zeta}}$$
such that 
$$\langle x,y\rangle=\re\ h(s(x),s(y)) \ \ \textrm{ for all } x,y\in C_1$$
with $\langle \cdot,\cdot\rangle$ and $h(\cdot,\cdot)$ as above.
\end{proposition}
\begin{proof}
The construction of $s$ was explained above, the remaining issue is to argue that it does not depend on the choices made when defining the lifts. There are two things that we have to check: the compatibility of (\ref{req2}) around a triangle $t_{ijk}$ and the compatibility with loops around a vertex $v_i$. Consider first a triangle $t_{ijk}$ oriented such that $e_{ij}$ is followed by $e_{jk}$ and then $e_{ki}$. Suppose that $\widetilde{e_{ij}}$ is given. Using lemma \ref{lemmal}, we define the lifts $\widetilde{e_{jk}}$ and $\widetilde{e_{ki}}$. We have to check that $\widetilde{e_{ij}}$ and $\widetilde{e_{ki}}$ satisfy (\ref{req2}). The lifts $\widehat{e_{ij}}$,  $\widehat{e_{jk}}$ and $\widehat{e_{ki}}$ used in \ref{lemmal} to define our lifts have the property that they form an oriented triangle on $C$, or more precisely $\widehat{e_{ij}}+\widehat{e_{jk}}+\widehat{e_{ki}}$ belongs to the image of $C_2'(\calE)_{\overline{\zeta}}$ in $C_1'(\calE)_{\overline{\zeta}}$ (by abuse of notation, we denote this image by $C_2'(\calE)_{\overline{\zeta}}$). Taking into account the scaling (\ref{scale}), it follows that
\begin{equation}\label{radicale}
\widetilde{e_{ij}}+\omega\cdot\widetilde{e_{jk}}+\omega^2\cdot\widetilde{e_{ki}}\in C_2'(\calE)_{\overline{\zeta}}.
\end{equation}
Since $h$ is $0$ on $C_2'(\calE)_{\overline{\zeta}}$, we get: 
$$h(\widetilde{e_{ij}},\widetilde{e_{ij}}+\omega\cdot\widetilde{e_{jk}}+\omega^2\cdot\widetilde{e_{ki}})=0.$$
Solving for  $h(\widetilde{e_{ij}}, \widetilde{e_{ki}})$ gives the expected value, fact equivalent to the  identity:
$$-(\cot x+\cot y)+\frac{e^{-x i}}{\sin x}+\frac{e^{y i}}{\sin y}=0.$$

Similarly, consider a cycle of edges $e_{ij_0},e_{ij_1},\dots,e_{ij_d}$ around a vertex $v_i$, with $j_d=j_0$ and such that the pairs $(e_{j_0i}, e_{ij_1})$, $(e_{j_1i},e_{ij_2})$, etc. are consecutive edges  in the respective triangles (w.r.t.  the orientation of $T$). Once $\widetilde{e_{ij_0}}$ is defined, we can use lemma \ref{lemmal} to define $\widetilde{e_{ij_1}}$. Next, to define $\widetilde{e_{ij_2}}$, the are two possibilities: to use either $\widetilde{e_{ij_0}}$ or $\widetilde{e_{ij_1}}$. We note first that the two are equivalent. Namely, since there is no loop around $v_i$, the lifts $\widehat{e_{ij_0}}$, $\widehat{e_{ij_1}}$ and $\widehat{e_{ij_2}}$ are defined unambiguously. We then clearly have (either way we use to define $\widetilde{e_{ij_2}}$):
\begin{eqnarray*}
\widetilde{e_{ij_0}}&=& \widehat{e_{ij_0}}\\
\widetilde{e_{ij_1}}&=& e^{\alpha_1i}\cdot  \widehat{e_{ij_1}}\\
\widetilde{e_{ij_2}}&=& e^{(\alpha_1+\alpha_2)i}\cdot  \widehat{e_{ij_2}}
\end{eqnarray*}
[N.B. $e^{-(\pi-\alpha)i}$ of (\ref{scale}) becomes $e^{\alpha i}$ due to the change of orientation $e_{ij}$ vs. $e_{ji}$]. Say for simplicity $d=5$. When we complete a full loop around $v_i$, by the above discussion we have
\begin{eqnarray*}
\widetilde{e_{ij_0}}&=& \widehat{e_{ij_0}}\\
\widetilde{e_{ij_d}}&=& e^{\frac{5\pi}{3}i}\cdot  \widehat{e_{ij_d}},
\end{eqnarray*}
but now $\widehat{e_{ij_0}}$ and $\widehat{e_{ij_d}}$ both come from $e_{ij_0}$, but from lifts differing by a sheet. Thus,  since with our orientation convention the arc $(e_{j_0j_1},e_{j_1j_2},\dots)$ is wrongly oriented, we get
$$\widehat{e_{ij_d}}=\zeta\cdot  \widehat{e_{ij_0}}$$ 
(the two lifts differ by $r^{-1}$ and they are $\overline{\zeta}$-eigencycles). We conclude
$$\widetilde{e_{ij_0}}=\overline{\zeta}\cdot \widehat{e_{ij_d}}=\overline{\zeta}\cdot \zeta \cdot \widehat{e_{ij_0}}=\widetilde{e_{ij_0}}$$
as needed.
\end{proof}

We recall that the lattice $R^\perp_{\overline{L}}$ has a presentation as the quotient of $C_1$ by $K$, where $K$ is the image of a map $C_0\to C_1$ sending the generators of $C_0$ to the elements $\zeta_i\in C_1$ (see (\ref{eqzeta})). We now note that we can write:
\begin{eqnarray*}
-\zeta_i&=&-(e_{ij_0}+e_{ij_1}+\dots+e_{ij_d})+(e_{j_0j_1}+e_{j_1j_2}+\dots+e_{j_{d-1}j_0})\\
&=&\left((\om+\omb)e_{ij_0}+\dots+(\om+\omb)e_{ij_{d-1}}\right)+(e_{j_0j_1}+e_{j_1j_2}+\dots+e_{j_{d-1}j_0})\\
&=&(e_{j_0j_1}+\om e_{ij_1}+\omb e_{ij_0})+\dots +(e_{j_{d-1}j_0}+\om e_{ij_0}+\omb e_{ij_{d-1}}).
\end{eqnarray*}
Using this decomposition (which is unique, when taking into account an orientation) and (\ref{radicale}),  we obtain (at least over $\bQ$): 
\begin{equation}\label{relation}
\xymatrix{
0\ar@{->}[r]&\bQ\ar@{->}[r]&C_0(\bQ)\ar@{->}[r]\ar@{->}[d]&C_1(\bQ)\ar@{->}[r]\ar@{->}[d]^{s}&R^\perp_{\overline{L}}\otimes \bQ\ar@{^{(}-->}[d]^{\bar{s}}\ar@{->}[r]&0\\
&0\ar@{->}[r]&C_2'(\bQ(i\sqrt{3}))_{\overline{\zeta}}\ar@{->}[r]&C_1'(\bQ(i\sqrt{3}))_{\overline{\zeta}}\ar@{->}[r]&\sH_1(C,\bQ(i\sqrt{3}))_{\overline{\zeta}}\ar@{->}[r]&0
}
\end{equation}
This induces an embedding $\bar{s}$ of the non-vanishing lattice $R^\perp_{\overline{L}}\subseteq \overline{L}$ that we associated to a triangulation via degenerate $K3$ surfaces into the lattice of cycles $\sH_1(C,\calE)_{\overline{\zeta}}$. The embedding is only an embedding of $\bZ$-lattices, i.e. we forget the complex structure of $\sH_1(C,\calE)_{\overline{\zeta}}$ and consider the real part of the Hermitian form. However, since
$$\widetilde{e_{ij}}+\omega\cdot\widetilde{e_{jk}}+\omega^2\cdot\widetilde{e_{ki}}\in C_2'(\calE)_{\overline{\zeta}},$$
when we consider the classes in $\sH_1(C,\calE)_{\overline{\zeta}}$ we obtain:
\begin{equation}\label{compatibil}
[\widetilde{e_{ij}}-\widetilde{e_{ki}}]=\omega\cdot [\widetilde{e_{ki}}-\widetilde{e_{jk}}].
\end{equation}
This can be interpreted as saying
$$
[\widetilde{\rho(e_{ki}-e_{jk})}]=[\widetilde{e_{ij}-e_{jk}}]=\omega\cdot [\widetilde{e_{ki}-e_{jk}}],$$
where $\rho$ is the order $3$ automorphism introduced in \S\ref{sectmu3} (see esp. Prop. \ref{mu3onc}). In other words, when we restrict $\bar{s}$ to the primitive part $Q$, the embedding $\bar{s}$ becomes an embedding of Eisenstein lattices, where the  Eisenstein lattice structure on $Q$ is given by construction of section \ref{secteisenstein}, and that on $H_1(C,\calE)_{\overline{\zeta}}$ is the natural structure [N.B. $Q=R^\perp_{\overline{L}}\cap P$ in the notation of \S\ref{sectmu3}].

\smallskip

Finally, we note that it is easy to adapt the discussion of this section to the case that we have hexagonal points. Namely, instead of geodesics joining two singular points of the triangulation, we have to lift pseudo-geodesics joining two singular points of the triangulation. This can be done in an obvious way. Namely, for the edges incident to singular points we proceed  as before. What remains are straight segments passing through the hexagonal points.  By the inductive procedure that we used earlier, we can assume that one end-point of the straight segment  was already lifted. We lift the straight segment to the corresponding sheet of the cover $C\to \bP^1$. One should also note that our hexagonal relations (\ref{eqhex}) are compatible with the relations coming from eigenchains (see diagram (\ref{eigen2})) [N.B. each hexagonal point gives by (\ref{take1}) an element of $C_0'(\calE)_{\overline{\zeta}}$, viewed over $\bR$, produces two linear relations].

\begin{corollary}\label{req3}
Let $\gamma_{ij}$ and $\gamma_{kl}$ be two pseudo-geodesics as constructed in \S\ref{sectpseudo} and $\widetilde{\gamma}_{ij}$ and $\widetilde{\gamma}_{kl}$ their lifts to $\sH_1(C,\calE)_{\overline{\zeta}}$. Then,
$$\langle \gamma_{ij}, \gamma_{kl}\rangle=  \re\ h(\widetilde{\gamma}_{ij}, \widetilde{\gamma}_{kl})=-\frac{4}{\sqrt{3}}\cdot \re\ H(\widetilde{\gamma}_{ij}, \widetilde{\gamma}_{kl})$$
where $\langle\cdot, \cdot \rangle$ is the intersection pairing on $\overline{L}$ and $H$ is the Hermitian form used by Looijenga \cite{loijl} for $\sH_1(C,\calE)_{\overline{\zeta}}$. Furthermore, if $\gamma$ is a linear combination of pseudo-geodesics such that $\langle \gamma,h\rangle=0$, then 
$$\omega\cdot \widetilde{\gamma}=\widetilde{\rho(\gamma)},$$
where $\rho$ denotes the order $3$ isometry constructed in section \ref{secteisenstein}. 
\qed
\end{corollary}

\subsection{Proof of Thm. \ref{mainthm2}}\label{proof2}
Assume for simplicity that there are only vertices of degree $5$ and $6$. Otherwise one has to replace $\overline{L}$ by the non-vanishing lattice $R^\perp_{\overline{L}}$ and to  keep track of the discriminant groups. Since the singularity behavior on the two sides is essentially identical, this can be done quite easily.

By the results of section \ref{sectpseudo}, the lattice $\overline{L}$ is generated by pseudo-geodesics (even over $\bZ$ with our assumption). In section \ref{sectlift}, we have constructed a $\bZ$-linear map 
\begin{equation}\label{embed1}
s:\overline{L}\hookrightarrow \sH_1(C,\calE)_{\overline{\zeta}}\cong M^\calE
\end{equation}
by lifting the pseudo-geodesics $\gamma_{ij}$ to eigencycles $\widetilde{\gamma_{ij}}$. We noted that when we forget the Eisenstein structure and take the associated bilinear form (see Sect. \ref{defeisenstein}), we obtain an isometric lattice embedding (see Cor. \ref{req3}). We also showed that when restricted to pseudo-geodesics of length  $0$ (i.e. $\gamma.h=0$), the Eisenstein structure introduced in section \ref{secteisenstein} coincides with the Eisenstein structure induced from $M^\calE$. 

\smallskip

Let $\delta_2^\calE\in M^\calE$ be the cocycle associated by Thurston (cf. Thm. \ref{thmthurston}). As noted by Looijenga \cite[Remark 1.1]{loijl}, the cocycle $\delta_2^\calE\in M^\calE$ can be represented by the (multi-valued) Lauricella differential $\eta$. Its pull-back $\widetilde{\eta}$ to $C$ is a single-valued differential form with class generating $\sH^{1,0}(C)_{\overline{\zeta}}$. Obviously, we can change $\eta$ by scaling. We claim that we can chose a representative for $\eta$ such that  for every $\bR$-linear combination $\gamma$ of pseudo-geodesics we have
\begin{equation}\label{equiv1}
\frac{1}{2}\cdot \langle \gamma,h\rangle=\frac{1}{6}\cdot \int_{\widetilde{\gamma}} \widetilde{\eta}= (\textrm{combinatorial length of } \gamma),
\end{equation}
where the left hand side is with respect to the intersection form on $\overline{L}$ and the right hand side is the usual evaluation of differential forms on cycles on $C$ (the scaling factor of $\frac{1}{6}$ is due to the cover). Considering a single pseudo-geodesic $\gamma_{ij}$, it is clear that  (\ref{equiv1}) can be achieved for some scaling of the Lauricella differential $\eta$. Namely, on one hand, $|\eta|^2$ gives locally the cone metric of Thurston (cf. \cite[Remark 1.1]{loijl}). On the other hand, by construction (see \S\ref{sectpseudo}), $\gamma_{ij}$ is locally a straight segment in the Euclidean plane. By scaling $\eta$ such that the length of an edge of the triangulation is $1$ and applying Lemma  \ref{lemmageo}, we get
\begin{equation}\label{equiv2}
\frac{1}{3}\cdot \frac{\int_{\widetilde{\gamma_{ij}}} \widetilde{\eta}}{\langle \gamma_{ij},h\rangle}\in S^1\subset \bC^*
\end{equation}
for each pseudo-geodesic $\gamma_{ij}$. By a rotation, we can assume that the above fraction is $1$ for some reference pseudo-geodesics $\gamma_{i_0j_0}$.

To prove that (\ref{equiv1}) holds always, by $\bR$-linearity, it suffices to show that 
\begin{equation}\label{equiv3}
\frac{1}{2}\left(\langle \gamma_{ij},h\rangle+ \langle \gamma_{jk},h\rangle\right)=\frac{1}{6}\left(\int_{\widetilde{\gamma_{ij}}} \widetilde{\eta}+\int_{\widetilde{\gamma_{jk}}} \widetilde{\eta}\right)
\end{equation}
for pairs of pseudo-geodesics meeting in a point. By an inductive procedure starting with $\gamma_{i_0j_0}$, we can assume 
$$\frac{1}{2}\cdot \langle \gamma_{ij},h\rangle=\frac{1}{6}\cdot \int_{\widetilde{\gamma_{ij}}}\widetilde{\eta}= (\textrm{length of } \gamma_{ij}).$$
We can also assume  (the other cases are similar) that locally at the intersection point $v_j$, the pseudo-geodesics $\gamma_{ij}$ and $\gamma_{jk}$ are consecutive edges of a standard equilateral triangle, oriented according to the counter-clockwise orientation. As in \S\ref{sectlift}, we have (locally at the intersection point):
\begin{eqnarray*}
\widetilde{\gamma_{ij}}&=& \epsilon'_j\\
\widetilde{\gamma_{jk}}&=& \omb\cdot \epsilon'_{j+1}
\end{eqnarray*}
where $\epsilon'_j$ and  $\epsilon'_{j+1}$ are lifts of consecutive arcs in a Lauricella arc system. Again referring to the metric interpretation of Thurston,  $\frac{1}{6}\cdot \int_{\epsilon'_{j+1}} \widetilde{\eta}\in \calE$ with norm equal to length of the pseudo-geodesic $\gamma_{jk}$, but with argument equal to the change of angle between $\epsilon'_j$ and  $\epsilon'_{j+1}$. We conclude 
\begin{equation}
\frac{1}{6}\cdot \int_{\widetilde{\gamma_{jk}}} \widetilde{\eta}=\frac{\omb}{6} \cdot \int_{\epsilon'_{j+1}} \widetilde{\eta}=\frac{\omb}{6} \cdot \om\cdot  \int_{\epsilon'_{j+1}} |\widetilde{\eta}|= (\textrm{length of } \gamma_{jk})
\end{equation}
as needed. In other words, the lifts of pseudo-geodesics are normalized such that the cocycle of Thurston is measuring their combinatorial length. The same is true about $\delta'=\frac{h-\delta}{2}$ (i.e. a generator of $M^\calE$ in our construction). Thus they must coincide. A priori, due to the scaling necessary to make the first integral real,  this is true only up to multiplying by elements of $\bQ(i\sqrt{3})\cap S^1$ (i.e. up to preserving the shape of the triangulation). However, since both $h$ and the cocyle of Thurston are constructed with the basic unit a standard equilateral triangle, a unit of $\calE$ suffices. Thus, indeed we get the same element in $M^\calE$, i.e.  both geometric constructions associate the same arithmetic invariant to a triangulation of non-negative combinatorial curvature. \qed

%%%%%%%%%%%%%%%%%%%%%%%%%%%%%%%%%%%%%%%%%%%%%%%%%%%%%%%%%%%%%%
\bibliography{references}
\end{document}